\documentclass[reqno]{amsart}
\usepackage{graphicx}
\usepackage{xcolor,bbm,easybmat,bm}
\usepackage{amsmath,amsfonts,amsthm,amscd,amssymb,mathrsfs,esint}
\usepackage{enumerate}
\usepackage[shortlabels]{enumitem}
\usepackage[colorlinks,citecolor=red,pagebackref,hypertexnames=false]{hyperref}

\DeclareMathOperator{\divg}{div}

\DeclareMathOperator{\osc}{osc\ }
\DeclareMathOperator*{\fiint}{\ensuremath{\iint\text{\kern-1.36em{\raisebox{5.87pt}{\rotatebox{-93}{$\setminus$}}}}}\,}

\theoremstyle{plain}
\newtheorem{theorem}[equation]{Theorem}

\newtheorem{lem}[equation]{Lemma}

\newtheorem{cor}[equation]{Corollary}

\theoremstyle{definition}

\newtheorem{defn}[equation]{Definition}

\theoremstyle{remark}

\newtheorem{re}[equation]{Remark}

\numberwithin{equation}{section}
\newcommand{\norm}[1]{\left\Vert#1\right\Vert}
\newcommand{\abs}[1]{\left\vert#1\right\vert}
\newcommand{\br}[1]{\left(#1\right)}
\newcommand{\set}[1]{\left\{#1\right\}}
\newcommand{\Real}{\mathbb R}
\newcommand{\Rn}{\mathbb R^n}

\newcommand{\eps}{\varepsilon}

\newcommand{\C}{\mathcal{C}}

\newcommand{\vp}{\varphi}

\newcommand{\wt}{\widetilde}
\newcommand{\bdy}{\partial}
\newcommand{\Dt}{\partial_t}
\newcommand{\Di}{\partial_i}
\newcommand{\Dj}{\partial_j}

\newcommand{\1}{\mathbbm{1}}

\newcommand{\RR}{{\mathbb{R}}}

\newcommand{\ZZ}{{\mathbb{Z}}}

\newcommand{\ep}{\hfill$\Box$ \vskip 0.08in}

\newcommand{\dv}{\operatorname{div}}

\newcommand{\R}{\mathbb{R}} 
\newcommand{\cN}{\mathfrak{N}} 
 
\newcommand{\cA}{\mathfrak{A}}
\newcommand{\sm}{\setminus} 
\renewcommand{\epsilon}{\varepsilon}

\renewcommand{\d}{\partial} 
\newcommand{\ms}{\medskip}

\begin{document}

\title
{Carleson measure estimates for the Green function}
\thanks{
G. David was partially supported by the European Community H2020 grant GHAIA 777822,
and the Simons Foundation grant 601941, GD.
S. Mayboroda was partly supported by the NSF RAISE-TAQS grant DMS-1839077 and the Simons foundation grant 563916, SM. 
}

\author{Guy David}
\author{Linhan Li}
\author{Svitlana Mayboroda}

\newcommand{\Addresses}{{
  \bigskip
  \footnotesize

 Guy David, \textsc{Universit\'e Paris-Saclay, CNRS, 
 Laboratoire de math\'ematiques d'Orsay, 91405 Orsay, France} 
 \par\nopagebreak
  \textit{E-mail address}: \texttt{guy.david@universite-paris-saclay.fr}

  \medskip

  Linhan Li, \textsc{School of Mathematics, University of Minnesota, Minneapolis, MN 55455, USA}\par\nopagebreak
  \textit{E-mail address}: \texttt{li001711@umn.edu}

  \medskip

  Svitlana Mayboroda, \textsc{School of Mathematics, University of Minnesota, Minneapolis, MN 55455, USA}\par\nopagebreak
  \textit{E-mail address}: \texttt{svitlana@math.umn.edu}

}}

\maketitle

\begin{abstract}
In the present paper we consider an elliptic divergence form operator in the half-space and 
prove that its Green function is almost affine, or more precisely, that the normalized difference between the Green function 
and a suitable affine function at every scale 
satisfies a Carleson measure estimate, provided that the oscillations of the coefficients satisfy the traditional quadratic Carleson condition. The results are sharp, and in particular, it is demonstrated that the class of the operators considered in the paper cannot be improved.  
\end{abstract}

\tableofcontents

\section{Introduction}
Let $L=-\divg(A\nabla)$ be a divergence form elliptic operator on the upper half-space 
$\Real^{d+1}_+$. In the present paper we show that if $L$ is reasonably well-behaved
then the Green function for $L$ is well approximated by multiples
of the distance to $\R^d$. There are many predecessors of these results which we will discuss below (\cite{kenig2001dirichlet},\cite{dindos2007lp},\cite{hofmann2017implies},\cite{hofmann2017uniform}  to mention only the closer ones).
At this point, however, let us underline two important points. First, the class of the operators that we consider is of the nature of the best possible, as shown by the counterexamples in Section~\ref{sec optm}. The estimates themselves are sharp, and in fact, a weak version of them is equivalent to the uniform rectifiability \cite{DM2020}. We hope to ultimately show that 
the much stronger estimate proved here is also true for domains with a uniformly rectifiable boundary,
thus giving a strong and a weak characterization of uniform rectifiability in terms of approximation of the
Green function (or more generally solutions) by distance function. But this will have to be the subject of another paper.
Secondly, the method of the proof itself is quite unusual for this kind of bounds. A typical approach is through 
integrations by parts,
which, however, does not allow one to access the optimal class of the coefficients.  Roughly speaking, we are working with the square of the second derivatives of the Green function and given the roughness of the coefficients, there are too many derivatives
in to control
to take advantage of the equation while integrating by parts. Here, instead, we make intricate comparisons with solutions of the constant-coefficient operators, carefully adjusting them from scale to scale. We feel that the method itself is a novelty for this circle of questions and that it illuminates the nature of the Carleson estimates in a completely different way, hopefully opening a door to many other problems.

More generally, we are interested in the relations between an elliptic operator $L$ 
on a domain $\Omega$, the geometry of $\Omega$, and the boundary behavior 
of the Green function. It is easy to see that the Green function with a pole at infinity for the Laplacian on the upper half-space $\Real^{d+1}_+:=\set{(x,t): x\in\Real^d, t\in \Real_+}$ 
is a multiple of $t$, the distance to the boundary, and more generally the Green 
function with a pole that is relatively far away is close to the distance function. There have been many efforts to generalize this to more general settings.  For instance, in \cite{caffarelli1981existence} the authors obtain flatness of the boundary from local small oscillations of the gradient of the Green function with a pole sufficiently far away. 
Philosophically,  similar considerations underpin the celebrated results of Kenig and Toro 
connecting the flatness of the boundary to the property that the logarithm of the Poisson kernel 
lies in VMO \cite{kenig1999free}. Much more close to our setting is the study of the so-called Dahlberg-Kenig-Pipher operators (defined in \eqref{1a7}-\eqref{1a8}) pioneered by Kenig and Pipher \cite{kenig2001dirichlet},\cite{dindos2007lp} in combination with the study of the harmonic measure on uniformly rectifiable sets by Hofmann, Martell, Toro, Tolsa, and others (see \cite{hofmann2017harmonic},\cite{hofmann2017uniform} and many of their predecessors). Undoubtedly, 
the behavior of the harmonic measure is connected to  the regularity of  Green function $G$, yet the latter is different and surprisingly has been much less studied. In part, this is due to the fact that the harmonic measure is related to the gradient of $G$ at the boundary while the estimates we target in this paper reach out to the second derivatives of $G$. One could say that the two are related by 
an integration by parts,
but in the world of the rough coefficients this is not so. Indeed, relying on these ideas, \cite{hofmann2017implies} establishes 
second derivatives estimates 
for the Green function somewhat similar to ours under a much stronger condition that the gradient of the coefficients, rather than its square, satisfies a Carleson condition. It was clear already then that the optimal condition must be a control of the square-Carleson norm, but their methods, using the aforementioned integration by parts,  did not give a possibility to overcome this restriction. In this paper we achieve the optimal results and, indeed, demonstrate using the counterexamples that they are the best possible.  

In the present paper, we focus on $\Omega = \Real^{d+1}_+$, and show that
for the operators satisfying a slightly weaker version of the Dahlberg-Kenig-Pipher condition 
described below, the Green function is well approximated by multiples of $t$, in the sense that
the gradient of normalized differences satisfies a square Carleson measure estimate.
Notice that the class of coefficients authorized below is enough to treat the case 
when $\Omega$ is a Lipschitz graph domain, by a change of variables. 
As we mentioned above, we plan  to pursue more general uniformly rectifiable sets in the upcoming work, which would give a much stronger version of our previous results in \cite{DM2020} and would show that our estimates are {\it equivalent} to the uniform rectifiability
of the boundary.
At this point, restricting to
the simple domain $\Omega = \Real^{d+1}_+$ will have the advantage of making the geometry
cleaner and focusing on one of the tools of this paper, concerning the dependence
of $G$ (or the solutions) on the coefficients. Even in the ``simple" case of the half-space,
the question of good approximation of $G$ by multiples of $t$ seems, to our surprise,
to be widely open, and the traditional methods of analysis break down 
brutally when trying to achieve such results. Perhaps one could also say that this setting is more classical. Let us pass to the details.

\ms
Consider an operator in divergence form $L=-\divg(A\nabla)$, 
where $A=\begin{bmatrix} a_{ij}(X) \end{bmatrix}$ is an $(d+1)\times (d+1)$ matrix of real-valued, bounded and measurable functions on $\Real^{d+1}_+$.
We say that $L$ is elliptic if there is some $\mu_0>1$ such that 
\begin{equation}
\langle A(X)\xi,\zeta\rangle \le \mu_0\abs{\xi}\abs{\zeta} \mbox{ and } \langle A(X)\xi,\xi\rangle \ge \mu_0^{-1}\abs{\xi}^2
\text{ for  }X \in \Real^{d+1}_+ \text{ and  }\xi, \eta\in\Real^{d+1}.
\label{cond ellp}
\end{equation}

We use lower case letters for points in $\Real^d$, e.g. $x\in\Real^d$, and capital letters 
for points in $\Real^{d+1}$, e.g. $X=(x,t)\in \Real^{d+1}$. 
We identify $\R^d$ with $\R^d \times \{ 0 \} \subset \R^{d+1}$ so, 
when $t=0$, we may write $x$ instead of $(x,0)\in\Real^{d+1}$.

For $x\in \R^d$ and $r > 0$, we denote by $\Delta(x,r)$ the surface ball
$B_r(x) \cap \set{t=0} \subset \R^d$. Thus $\Delta(x,r)$ is a ball in $\R^d$
while $B(x,r)$ is the ball of radius $r$ in $\R^{d+1}$.
We denote by  
\begin{equation} \label{TT}
T(x,r):=B_r(x)\cap \Real^{d+1}_+
\quad\text{and} \quad W(x,r):=\Delta(x,r)\times\Bigl(\frac{r}{2},r\Bigr] \subset \Real^{d+1}_+
\end{equation}
the corresponding Carleson box and Whitney cube.
Note that $T(x,r)$ is a half ball in $\Real^{d+1}_+$ over $\Delta(x,r)$. 
We may simply write $T_\Delta$ for a half ball over $\Delta\subset\Real^d$.

\begin{defn}[Carleson measure]\label{d13}
We say that a nonnegative Borel measure $\mu$ is a Carleson measure in $\Real_+^{d+1}$, 
if its Carleson norm 
\[\norm{\mu}_{\mathcal{C}}:=\underset{\Delta\subset \Real^d}{\sup}\frac{\mu(T_\Delta)}{\abs{\Delta}}\]
is finite, where 
the supremum is over all the surface balls $\Delta$ and 
$\abs{\Delta}$ is the Lebesgue measure of $\Delta$ in $\Real^d$. We use $\mathcal{C}$ to denote the set of Carleson measures on $\Real^{d+1}_+$. 

For any 
surface ball
$\Delta_0\subset\Real^d$, we use $\mathcal{C}(\Delta_0)$ to denote the set of Borel measures satisfying the Carleson condition restricted to $\Delta_0$, 
i.e., such that
\[\norm{\mu}_{\mathcal{C}(\Delta_0)}:=\underset{\Delta\subset \Delta_0}{\sup}\frac{\mu(T_\Delta)}{\abs{\Delta}} < +\infty
.\]
\end{defn}

Next we want to define a (weaker) version of the Dahlberg-Kenig-Pipher conditions in the form which is convenient for the point of view taken in this paper. 
We would like to say that the matrix $A = A(X)$ is often
close to a constant coefficient matrix. The simplest way to measure this is to use the numbers
\begin{equation} \label{1a4}
\alpha_\infty(x,r) = \inf_{A_0 \in \cA_0(\mu_0)} \, \sup_{(y,s) \in W(x,r)} |A(y,s)-A_0|,
\end{equation}
where the infimum is taken over the class $\cA_0(\mu_0)$ of (constant!) 
matrices $A_0$ that satisfy the ellipticity condition \eqref{cond ellp}. 
Notice that the matrix $A_0$ is allowed
to depend on $(x,r)$, so $\alpha_\infty(x,r)$ is a measure of the oscillation of
$A$ in $W(x,r)$, similarly to \cite{dindos2007lp}. We require $A_0$ to satisfy \eqref{cond ellp} for convenience,
but if we did not, we could easily replace $A_0$ by one of the $A(y,s)$,
$(y,s)\in W(x,r)$, which satisfies \eqref{cond ellp} by definition, at the price of multiplying 
$\alpha_\infty(x,r)$ by at most $2$. The same remark is valid for the slightly more general 
numbers 
\begin{equation} \label{1a5}
\alpha_q(x,r) = \inf_{A_0 \in \cA_0(\mu_0)}
\bigg\{\fint_{(y,s) \in W(x,r)} |A(y,s)-A_0|^q \bigg\}^{1/q}
\end{equation}
where in fact $q$ will be chosen equal to 2.

\begin{defn}[Weak DKP condition] \label{d1a6}
We say that the coefficient matrix $A$ satisfies the weak DKP condition with
constant $M > 0$, when $\alpha_2(x,r)^2 \frac{dxdr}{r}$ is a 
Carleson measure on $\R^{d+1}_+$, with norm
\begin{equation}\label{1a7}
\cN_2(A) : = \norm{\alpha_2(x,r)^2 \frac{dxdr}{r}}_{\mathcal{C}} \le M.
\end{equation}
\end{defn}

We may also say that $\alpha_2(x,r)^2$ satisfies a Carleson measure estimate. 
Recall that this implies that $\alpha_2(x,r)^2$ is small most of the time
(to the point of being integrable against the infinite invariant measure $\frac{dxdr}{r}$),
but does not vanish at any specific speed given in advance. 

The name comes from a condition introduced by Dahlberg, Kenig, and Pipher, which
instead demands that $\wt\alpha(x,r)^2$ satisfy a Carleson estimate, where
\begin{equation}\label{1a8}
\wt\alpha(x,r) = r \sup_{(y,s) \in W(x,r)} |\nabla A(y,s)|.
\end{equation}
In 1984, Dahlberg first introduced this condition, and conjectured  that such a Carleson condition guarantees 
the absolute continuity of the elliptic measure
with respect to the Lebesgue measure in the upper half-space. In 2001, Kenig and Pipher \cite{kenig2001dirichlet} proved Dahlberg's conjecture. Since it is obvious that $\alpha_2(x,r)\le \alpha_\infty(x,r) \leq 2\wt\alpha(x,r)$, 
we see that our condition is weaker than the classical 
DKP condition, but importantly they have the same homogeneity. A similar weakening of the DKP condition, pertaining to  the oscillations of the coefficients, has been considered, e.g. in \cite{dindos2007lp}.
We could also have chosen an exponent $q\in(2,\infty]$ for $\alpha_q$ in Definition \ref{d1a6}, but there is no point doing so as the H\"older inequality implies that the current condition is the weakest. Surprisingly, our theorem is easier to prove under this weaker condition. 

\ms
We now say what we mean by good approximation by affine functions. On domains other
than $\R^{d+1}_+$, we would use other models than the function $(y,t) \mapsto t$,
such as (functions of) the distance to the boundary, but here we are interested in 
(approximation by) the affine function $(y,t) \mapsto \lambda t$, with $\lambda > 0$.

We said earlier that we wanted to study the approximation of the Green functions 
(and we did not mention the poles too explicitly), but in fact our properties will also be valid
for positive solutions $u$ of $Lu = 0$ that vanish at the boundary. 

In addition, given such a solution $u$, when we are considering a given Carleson box $T(x,r)$, 
we do not want to assume any a priori knowledge on the 
average size of $u$ in $T(x,r)$,
so we just want to measure
the approximation of $u$, in $T(x,r)$, by the best affine function $a_{x,r}$ that we can think of,
and it is reasonable to pick
\begin{equation} \label{1a9}
a_{x,r}(z,t) = \lambda_{x,r} t, \, \text{ where } 
\lambda_{x,r} = \lambda_{x,r}(u) = \fint_{T(x,r)}\partial_t u(z,t)dzdt
\end{equation}
is the average on $T(x,r)$ of the vertical derivative. 
See the beginning of Section~\ref{sec const} for more details about this choice of $\lambda_{x,r}$. 
We measure the proximity of the two functions by the $L^2$ average of the difference 
of the gradients (we seem to forget $u$ but after all, it is easy to recuperate the functions 
from their gradients because they both vanish on the boundary), which we divide 
by the local energy of $u$ because we do want the same result
for $u$ as for $\lambda u$. That is, we set
\begin{multline}\label{1a10}
J_u(x,r) =  \fint_{T(x,r)} |\nabla_{z,t} (u(z,t)-a_{x,r}(z,t))|^2dzdt\\
= \fint_{T(x,r)} |\nabla_{z,t} u(z,t) - \lambda_{x,r}(u) \mathbf{e}_{d+1}|^2dzdt,
\end{multline} 
where $\mathbf{e}_{d+1} = (0,\ldots, 1)$ is the vertical unit vector, and then 
divide by 
\begin{equation}\label{def E}
    E_u(x,r)=\fint_{T(x,r)} |\nabla u|^2
\end{equation}
to get the number
\begin{equation} \label{1a11}
\beta_u(x,r) = \frac{J_u(x,r)}{E_u(x,r)}.
\end{equation}
This number measures the normalized non-affine part of the energy of $u$
in $T(x,r)$. We want to say that $u$ is often close $a_{x,r}$, i.e., that
$\beta_u(x,r)$ is often small, and this will be quantified
by a Carleson measure condition on $\beta_u$. We won't need to square $\beta_u$,
because $J_u$ is already quadratic.

The simplest version of our main result is the following.

\begin{theorem}\label{mt1}
Let $A$ be a $(d+1)\times (d+1)$ matrix of real-valued functions on 
$\Real^{d+1}_+$ 
satisfying the ellipticity condition \eqref{cond ellp}. 
If $A$ satisfies the weak DKP condition with some constant $M\in(0,\infty)$,
and if we are given $x_0 \in \R^d$, $R>0$, and a positive solution $u$ of 
$Lu=-\divg\br{A\nabla u}=0$ in $T(x_0,R)$, with $u=0$ on $\Delta(x_0,R)$, 
then the function $\beta_u$ defined by \eqref{1a11} satisfies a Carleson
condition in $T(x_0,R/2)$, and more precisely
\begin{equation}\label{1a13}
\norm{\beta_u(x,r) \frac{dxdr}{r}}_{\mathcal{C}(\Delta(x_0,R/2))} \leq C+C\,M
\end{equation}
where $C$ depends only on $d$ and  $\mu_0$.
\end{theorem}

That is, $u$ is locally well approximated by affine functions in $T(x_0,R/2)$,
with essentially uniform Carleson bounds. Here ``solution'' means ``weak solution'', and the values of $u$ on $\R^d$
are well defined because solutions are locally H\"older continuous up to the boundary; 
this will be explained better in the next section.

Notice that the constant $M>0$ can take any values, and we explicitly underlined the norm dependence. The result applies when $u$ is the Green function for $L$, 
with a pole anywhere in $\R^{d+1}_+ \sm \overline T(x_0,R)$. Even in the case 
of the Laplacian, the smallness of $M$ does not guarantee the smallness of \eqref{1a13}, that is, $u$ is not necessarily so close to an affine function at the scale $R$. 
This is natural (the impact  of what happens outside of $T(x_0,R)$ could be substantial), and this effect will be ameliorated in the next statement, at the price of some additional
quantifiers; the point is that the Green function with a pole at $\infty$, or even a positive
solution in a much larger box than $T(x_0,R)$, behaves better and has a better approximation. 
The next theorem says that we can have Carleson norms for $\beta_u$
that are as small as we want, provided that we take a small DKP constant and 
a large security box where $u$ is a positive solution that vanishes on the boundary.

 \begin{theorem}\label{mt2}
Let $d$, $\mu_0$ be given, let $u$ and $\Delta(x_0,R)$ satisfy the assumptions of Theorem \ref{mt1}, 
and let $A$ satisfy the weak DKP condition in $\Delta(x_0,R)$.
Then for $\tau \leq 1/2$ we have the more precise estimate
\begin{equation}\label{1a15} 
\norm{\beta_u(x,r) \frac{dxdr}{r}}_{\mathcal{C}(\Delta(x_0, \tau R))} 
\leq C \tau^a + C \norm{\alpha_2(x,r)^2 \frac{dxdr}{r}}_{\mathcal{C}(\Delta(x_0,R))},
\end{equation}
where $C$ and $a > 0$ depends only on $d$ and $\mu_0$.
\end{theorem}

This way the right-hand side can be made as small as we want.
Notice that we only need $A$ to satisfy
the weak DKP condition in $\Delta(x_0,R)$; the values of $A$ outside of
$T(x_0,R)$ should be irrelevant anyway, because we do not know anything about $u$ there.

We observed earlier that this result applies to the Green function with a pole at $\infty$ (see Lemma \ref{lem emGreen_infinity} for the precise definition), and to operators that satisfy the classical 
Dahlberg-Kenig-Pipher condition where the square of the function $\wt \alpha$ of \eqref{1a8}
satisfies a Carleson measure estimate. Notice that when $u$ is the Green function with pole at $\infty$ for $L$, Theorem \ref{mt2} implies that the Carleson norm of $\beta$ is simply less than $C \cN_2(A)$, with $\cN_2(A)$ as in \eqref{1a7}.

A rather  direct consequence of our results is a Carleson measure estimate on the second derivatives 
of the Green function for DKP operators. 
\begin{cor}\label{cor main}
Let $A$ be a $(d+1)\times (d+1)$ matrix of real-valued functions on 
$\Real^{d+1}_+$ satisfying the ellipticity condition \eqref{cond ellp}. 
Suppose $A$ satisfies the classical DKP condition with constant $C_0\in(0,\infty)$, that is, 
\begin{equation}\label{eq DKP}
    \norm{\wt{\alpha}(x,r)^2\frac{dxdr}{r}}_{\mathcal{C}}\le C_0,
\end{equation}
 where $\wt\alpha(x,r)$ is defined in \eqref{1a8}. If we are given $x_0 \in \R^d$, $R>0$, and a positive solution $u$ of 
$Lu=-\divg\br{A\nabla u}=0$ in $T(x_0,R)$, with $u=0$ on $\Delta(x_0,R)$, 
then there exists some constant $C$ depending only on $d$, $\mu_0$ and $C_0$ such that
\begin{equation}
    \int_{T_\Delta}\frac{\abs{\nabla^2u(y,t)}^2}{u(y,t)^2}\,t^3\,dydt\le C\abs{\Delta}
\end{equation}
for any $\Delta\subset\Delta(x_0,R/2)$.
\end{cor}
 We state this corollary on the upper half-space for simplicity, but it can be generalized to Lipschitz domains by a change of variables that preserves the DKP class operators. In fact, the change of variables will be a bi-Lipcshitz mapping whose second derivatives satisfy a Carleson measure estimate. With such regularity of the change of variables, as well as our estimates for $\beta_u$ in the main theorems, it reduces to the case of the upper half-space.

\ms
In Section \ref{sec optm}, we construct an operator that does not satisfy the DKP condition, for which the precise approximation estimates of Theorems~\ref{mt1} and \ref{mt2} fail.

In conclusion, let us point out that we extend the results above 
to domains with lower dimensional boundaries in \cite{david2021lower}.
In that case, there are currently no known free boundary results, in particular, 
it is not known whether the absolute continuity of elliptic measure with respect to the Hausdorff measure, or square function estimates, or the 
well-posedness of the Dirichlet problem imply the
rectifiability of the boundary, and we hope that the correct condition is, in fact, 
an analogue of the property that the Green function is almost affine. 
The first and the third authors of the paper started such a study in \cite{DM2020}, but if we want
precise approximation results for the Green functions, the first significant 
step in the positive direction should 
be a version of main results of the present paper in the higher co-dimensional context, 
and their extension to uniformly rectifiable sets.

The rest of this paper is organized as follows. In the next section we recall
some notation and the general properties of solutions that we need later.
In Section \ref{sec const} we comment the definition of $J_u$ and $\beta_u$,
prove some decay estimates for $\beta_u$ when $u$ is a weak solution of a constant 
coefficient operator, and extend this to the general case with a variational argument.
The rest of the proof of our main theorems, which consists in Carleson measure estimates with no special
relations with solutions, is done in Section \ref{sec cm}. We prove Corollary \ref{cor main} in Section \ref{sec cor} using Theorem \ref{mt1} and a Caccioppoli type argument. 
In Section \ref{sec optm}, we discuss the optimality of our results.

\section{Preliminaries 
and properties of the weak solutions
}\label{sec nota def}

In this section we recall some classical results for solutions of elliptic operators in divergence form.

Recall the notation $B(X,r)$ for open balls centered at $X\in \R^{d+1}$, $\Delta(x,r)$
for surface balls, $T(x,r)$ for Carleson boxes, and $W(x,r)$ for Whitney cubes (see near \eqref{TT}). Also denote by $\fint_Bf(x)dx:=\frac{1}{\abs{B}}\int_Bf(x)dx$ the average
of $f$ on a set $B$.

Let us collect some well-known estimates for solutions of $L=-\divg(A\nabla)$, where $A$ is a matrix of real-valued, measurable and bounded functions, satisfying the ellipticity condition \eqref{cond ellp}.

\begin{defn}[Weak solutions]\label{def weak sol}
Let $\Omega$ be a domain in $\Rn$. A function $u\in W^{1,2}(\Omega)$ is a weak solution to $Lu=0$ in $\Omega$ if for any $\vp\in W^{1,2}_0(\Omega)$, 
\[
\int_{\Omega}A(X)\nabla u(X)\cdot\nabla\vp(X)dX=0.
\]
\end{defn}

We will only be interested in the simple domains $\Omega = \R^{d+1}_+$ and  
$\Omega = \R^{d+1}_+ \cap B(x,r)$, with $x\in \R^d$ and $r>0$. The space $W^{1,2}_0(\Omega)$
is the closure in $W^{1,2}(\Omega)$ of the compactly supported smooth functions in $\Omega$.
Conventional or strong solutions are obviously weak solutions as well. 
In this paper, our solutions are always taken 
in the sense of Definition \ref{def weak sol}. 

From now on, $u$ is a (weak) solution in $\Omega$. When we say that $u=0$ on 
some surface ball $\Delta = \Delta(x,r) \subset \Omega$,
we mean this in the sense of $W^{1,2}(T_{\Delta})$. This 
means that $u$ is a limit in $W^{1,2}(T_{\Delta})$ of a sequence of functions in $C_0^1(\overline{T_{\Delta}}\setminus\Delta)$. 
We could also say that the trace of $u$, which is defined and lies in $H^{1/2}(\Delta)$,
is equal to $0$ on $\Delta$.
Ultimately, the De Giorgi-Nash-Moser theory (cf. Lemma \ref{lem bdy reg}) shows that under this assumption, the weak solution $u$ is in fact continuous in $T_{2r}\cup\Delta_{2r}$, and, in particular, $u$ vanishes on $\Delta$. 
Hence, in the rest of this paper the distinction is immaterial, but for now we will try to be precise.

We refer the readers to \cite{kenig1994harmonic} for proofs and references for the 
following lemmas.

\begin{lem}[Boundary Caccioppoli Inequality]\label{lem bdy cacio} 
Let $u\in W^{1,2}(T(x,2r))$ be a solution of $L$ in $T(x,2r)$, with $u=0$ on $\Delta(x,2r)$. 
There exists some constant $C$ depending only on the dimension and the ellipticity constant of $L$, such that
 \[
 \fint_{T(x,r)}\abs{\nabla u(X)}^2dX\le\frac{C}{r^2}\fint_{T(x,2r)}\abs{u(X)}^2dX.
 \] 
\end{lem}

\begin{lem}[Boundary De Giorgi-Nash-Moser inequalities]\label{lem bdy reg}
 Let $u$ be as in Lemma~\ref{lem bdy cacio}. Then 
 \[
 \sup_{T(x,r)}\abs{u}\le C \br{\fint_{T(x,2r)}u(X)^2dX}^{1/2}, 
 \]
 where $C=C(d,\mu_0)$. Moreover, for any $0<\rho<r$, we have, for some 
 $\alpha=\alpha(d,\mu_0)\in(0,1]$,
 \[
 \underset{T(x,\rho)}{\osc}u 
 \le C\br{\frac{\rho}{r}}^\alpha\br{\fint_{T(x,2r)}u(X)^2dX}^{1/2}, 
 \]
 where $\underset{\Omega}{\osc}u:=\underset{\Omega}{\sup\,}u-\underset{\Omega}{\inf\,}u$.
\end{lem}

\begin{lem}[Boundary Harnack Inequality]\label{lem bdy Harnack}
Let $u\in W^{1,2}(T(x,2r))$ be a nonnegative solution of $L$ in $T(x,2r)$ 
with $u=0$ on $\Delta(x,2r)$. 
Then 
\[
u(X)\le Cu(X_r) \qquad\forall\, X\in T(x,r), 
\]
where $C>0$ depends only on the dimension and $\mu_0$.
\end{lem}

Of course, each of these statements has an interior analogue where we would replace $T(x,r)$ 
by a ball $B(X,r)$ such that $B(X,2R) \subset \Omega$ and we 
would not have to specify the boundary conditions. The interior Harnack inequality reads 
as follows.
 
\begin{lem}[Harnack Inequality]
There is some constant $C$, depending only on the dimension and the ellipticity constant
for $A$, such that if $u\in W^{1,2}(\Omega)$ is a nonnegative solution of $Lu=0$ in 
$B(X,2r)\subset\Omega$, then
\[
\sup_{B(X,r)} u\le C\inf_{B(X,r)} u. 
\]
\end{lem}

We will also use the Comparison Principle.

\begin{lem}[Comparison Principle] 
Let $u,v\in W^{1,2}(T(x,2r))$ 
be two nonnegative solutions of $L$ in $T(x,2r)$, such that  
$u=v=0$ on $\Delta(x,2r)$ and $v$ is not identically null. 
Set $X_{x,r} = (x,r)$ (a corckscrew point for $T(x,2r)$). Then 
\[
C^{-1}\frac{u(X_{x,r})}{v(X_{x,r})} 
\le\frac{u(X)}{v(X)}\le C \frac{u(X_{x,r})}{v(X_{x,r})} 
\quad \text{ for all } X\in T(x,r), 
\]
where $C=C(n,\mu_0)\ge 1$.
\end{lem}

\begin{lem}[Reverse H\"older Inequality on the boundary]\label{lem RH}
We can find an exponent $p > 2$ and a constant $C \geq 1$, that depend only on 
$d$ and the ellipticity constant $\mu_0$ for $A$, such that if $u$ and $T(x,2r)$ are as 
in Lemma \ref{lem bdy cacio}, then 
\[
\br{\fint_{T(x,r)}\abs{\nabla u(X)}^pdX}^{1/p} 
\le C\br{\fint_{T(x,2r)}\abs{\nabla u(X)}^2dX}^{1/2}. 
\]
\end{lem}

See \cite{giaquinta1983multiple}, Chapter V for the proof of this Lemma. 

We prove the following simple consequence of the above for reader's convenience. 

\begin{lem}\label{lem corkscrew}
Let $u\in W^{1,2}(T(x,R))$ be a nonnegative solution of $L$ in $T(x,R)$,
with $u=0$ on $\Delta(x,R)$. 
Then for all $0<r<R/2$,
\begin{equation}\label{eqcs1} 
\fint_{T(x,r)}\abs{\nabla u(X)}^2dX 
\approx
\frac{u^2(X_{x,r})}{r^2}, 
\end{equation}
where $X_{x,r} = (x,r)$ as above and
the implicit constant depends only on $d$ and $\mu_0$.
\end{lem}

\begin{proof}
By translation invariance, we may assume that $x_0$ is the origin.

To prove the $\gtrsim$ inequality in \eqref{eqcs1}, we apply Lemma \ref{lem bdy reg},  
Lemma \ref{lem bdy Harnack}, and the Poincar\'e inequality, and get 
\[
u(X_{x,r})^2 \le C\,\sup_{T_{x, r/2}} 
u^2\le C\fint_{T_{x,r}}u^2(X)dX \le C r^2 \fint_{T_{x,r}} 
\abs{\nabla u}^2.
\]
For the $\lesssim$ inequality in \eqref{eqcs1}, simply 
combine the boundary Caccioppoli and boundary Harnack inequalities.
\end{proof}

\ms
We now record a basic regularity estimate for constant coefficient operators.
This will be used in the next section to get decay estimates for $J_u$, 
and then extended partially to our more general operators $L$, with comparison arguments. 
We shall systematically use $A_0$ to denote a {\it constant} real $(d+1)\times (d+1)$ matrix, 
which we always assume to satisfy the ellipticity condition \eqref{cond ellp},
and write $L_0:= -\divg\br{A_0\nabla}$. Solutions to such operators enjoy additional regularity and in particular, we will use the following result. 
We state it in $T_1 = T(0,1)$ to simplify the notation.
More generally, set $T_r = T(0,r)$ for $r > 0$.

\begin{lem}\label{l2a10}
Let $u\in W^{1,2}(T_1)$ be a solution to $L_0u=0$ in $T_1$ with $u=0$ on $\Delta_1$. 
Then for any multiindex $\alpha$, $\abs{\alpha}\in\mathbb{Z}_+$,
\begin{equation}\label{eq reg solL0}
 \underset{T_{\frac{1}{2}}}{\sup}\abs{D^{\alpha}u}
 \le C\br{\fint_{T_1}\abs{\nabla u(X)}^2dX}^{1/2},
\end{equation}
where $C=C(d,\mu_0,\abs{\alpha})$.
In particular, for any $T(x,r)\subset T_{1/2}$,
\begin{equation}\label{eq osc in use}
 \underset{T(x,r)}{\osc}\Di u\le Cr\br{\fint_{T_{1}} 
 \abs{\nabla u(X)}^2dX}^{1/2}, \quad i=1,2,\dots, d+1, 
\end{equation}
where the constant $C$ depends only on the dimension and $\mu_0$.
\end{lem}

\begin{proof} First we claim that 
the standard local estimates on solutions for constant-coefficient 
operators in $\RR^{d+1}_+$ ensure that 
\begin{equation}\label{eqloc}
\|D^\alpha u\|_{L^2 (T_{1/2})}\lesssim \|\nabla u\|_{L^2 (T_1)}+ \|u\|_{L^2 (T_1)}.
\end{equation}
This is due to the fact that any 
weak solution to $Lu=f$ on a smooth bounded domain $\Omega$ and 
with zero Dirichlet boundary data satisfies 
$$
\|u\|_{W^{m+2,2}(\Omega)} 
\lesssim \|f\|_{W^{m,2}(\Omega)} +\|u\|_{L^2(\Omega)}, \quad m=0,1,2,...; 
$$
see, e.g., \cite{evans2010partial}, \S~6.3, Theorems 4, 5.  
Here, $W^{m,2}(\Omega)$ is the Sobolev space of functions whose 
derivatives up to the order $m$ lie in $L^2(\Omega)$. With this at hand, we observe that 
for any smooth cutoff function $\eta$ equal to 1 on 
$B_{1/2}$ 
and supported in $B_{3/4}$ we have 
$$L_0(u\eta)=-A_0 \nabla \eta\cdot \nabla u-A_0\nabla u\cdot\nabla\eta +u\,L_0 \eta, $$
and hence the estimate above applied consecutively with $m=0,1,2...$ in some smooth domain 
$T_{3/4}\subset \Omega \subset T_1$ gives \eqref{eqloc}. Applying  
Poincar\'e's inequality, we conclude that 
\begin{equation}\label{eqloc2}\|D^\alpha u\|_{L^2 (T_{1/2})} 
\lesssim \|\nabla u\|_{L^2 (T_1)} 
\end{equation}
for any multiindex $\alpha$ with $|\alpha|\in \ZZ_+$. On the other hand, by the Sobolev embedding theorem (\cite{adams2003sobolev} Theorem 4.12), for any multiindex $\alpha$,
\[
\sup_{T_{1/2}}\abs{D^{\alpha}u}\le C \norm{u}_{W^{|\alpha|+n,2}(T_{1/2})} ,
\]
where $C$ depends on  $n$ and $|\alpha|$.
We combine this with \eqref{eqloc2} and get \eqref{eq reg solL0}.

The estimate \eqref{eq osc in use} is an immediate consequence of \eqref{eq reg solL0}, since  
\[
\underset{T(x,r)}{\osc} \partial_i u \le r \sup_{T(x,r)} 
\abs{\nabla \partial_i u}
\le r\,\sup_{T_{1/2}}\abs{\nabla \partial_i u}
\le Cr \br{\fint_{T_{1}}\abs{\nabla u}^2}^{1/2}, 
\]
as desired.
\end{proof}

\begin{re}\label{r2}
Lemma \ref{l2a10} is more than enough to prove Theorems \ref{mt1} and \ref{mt2}
in the special case of constant-coefficient operators. Indeed it says that $\nabla u$
is Lipschitz in $T_{1/2}$, so in particular $\nabla u - \nabla u(0)$ is small near the origin.
Notice that $\nabla u(0) = (0,\d_t u(0))$ because $u$ vanishes on the boundary;
with this and similar statements for other surface balls, it would be rather easy to control
$\beta_u$ and prove the theorems in the case of constant-coefficient operators. We don't do this here because we need more general estimates anyway.
\end{re}

\section{Approximations and the main conditional decay estimate}
\label{sec const}

We observed in Remark \ref{r2} that our theorems should be easy to prove when $L$
is a constant coefficient operator. In this section, we use the results of the previous section, 
together with an approximation argument, to prove some decay estimate for $\beta_u$
in regions where $A$ is nearly constant. See Corollary \ref{cor itr}.

At the center of the proof is an estimate for $||\nabla u - \nabla u_0||_2$, 
where $u$ is a solution for $L$ in some Carleson box $T(x,r)$, 
and $u_0$ is a solution for a close enough constant coefficient operator $L_0$, 
with the same boundary values on $\d T(x,r)$. See Lemma~\ref{lem comp u u0}.

\subsection{A little more about orthogonality, $J_u$, and $\beta_u$}\label{subsec orth}
First return to the approximation of a solution $u$ by the affine function 
$a_{x,r}(z,t) = \lambda_{x,r} t$ of \eqref{1a9}. Let us check what we said earlier,
that $a_{x,r}$ is the best affine approximation of this type in $T(x,r)$. Recall from 
\eqref{1a10} that
\begin{equation} \label{3a1}
\begin{aligned}
J_u(x,r) &=  \fint_{T(x,r)} |\nabla (u(z,t)-a_{x,r}(z,t))|^2dzdt
= \fint_{T(x,r)} |\nabla u - \lambda_{x,r}(u) \mathbf{e}_{d+1}|^2dzdt
\cr& = \fint_{T(x,r)}  |\nabla_z u(z,t)|^2dzdt + \fint_{T(x,r)} |\d_t u(z,t) - \lambda_{x,r}(u)|^2dzdt
\end{aligned}
\end{equation}
where $\mathbf{e}_{d+1} = (0,\ldots, 1)$ is the vertical unit vector, and we split the full 
gradient $\nabla u$ into the horizontal gradient $\nabla_x u$ and the vertical part $\d_t u$. 
Now $\lambda_{x,r}(u) = \fint_{T(x,r)} \d_t u$ by \eqref{1a9},
so $\d_t u - \lambda_{x,r}(u)$
is orthogonal to constants in $L^2(T(x,r))$, hence for any other $\lambda$, 
$$
\fint_{T(x,r)} |\d_t u - \lambda|^2 
= |\lambda - \lambda_{x,r}(u)|^2 + \fint_{T(x,r)} |\d_t u - \lambda_{x,r}(u)|^2,
$$ 
and, by the same computation as above,
\begin{equation} \label{3a2}
\begin{aligned}
\fint_{T(x,r)} |\nabla (u- \lambda t)|^2
&= |\lambda - \lambda_{x,r}(u)|^2 + \fint_{T(x,r)} |\nabla u - \lambda_{x,r}(u) \mathbf{e}_{d+1}|^2
\cr&= |\lambda - \lambda_{x,r}(u)|^2 + J_u(x,r).
\end{aligned}
\end{equation}
We may find it convenient to use the fact that, as a consequence,
\begin{equation} \label{3a3}
\beta_u(x,r) = \inf_{\lambda \in \R} \, 
\frac{\fint_{T(x,r)} |\nabla (u- \lambda t)|^2}{\fint_{T(x,r)} |\nabla u|^2} \leq 1.
\end{equation}
(compare with \eqref{1a11}, and for the second part try $\lambda = 0$).

\ms
For most of the rest of this section, we concentrate on balls centered at the origin; 
to save notation, we set $B_r = B(0,r)$, $T_r = T(0,r) = B_r \cap \R^{d+1}_+$, and 
$W_r = W(0,r)$ (see \eqref{TT}). Similarly, it will be convenient to use the notation
\begin{equation*}
J_u(r) = J_u(0,r)= \fint_{T_r}\abs{\nabla\br{u(x,t)-\lambda_r(u)\, t}}^2dxdt,
\end{equation*}
where \[\lambda_r(u) = \lambda_{0,r}(u) =\fint_{T_r} \partial_s u(y,s)dyds\]  
(see \eqref{1a9} and \eqref{1a10}). And we set \(E_u(r) = E_u(0,r)\),  \(\beta_u(r) = \beta_u(0,r)\) (see \eqref{def E} and \eqref{1a11}).

\subsection{Decay estimates for constant-coefficient operators}
We shall now prove a few estimates on solutions of constant-coefficient equation,
which will be useful when we try to replace $L$ by a constant-coefficient operator.
We start with a consequence of Lemma \ref{l2a10}. 

\begin{lem}\label{lem u0-lambda t}
Let $A_0$ be a constant matrix that satisfies the ellipticity condition \eqref{cond ellp},
set $L_0 = -\divg\br{A_0\nabla}$, and and let $u$ be a solution to $L_0 u = 0$ in $T_1$ 
such that $u=0$ on $\Delta_1$. 
There exists some constant $C$, 
depending only on the dimension and $\mu_0$, such that for
$0<r<1/2$,
\begin{equation}\label{eq u0-lambda t}
J_u(r)\le Cr^2J_u(1)\le Cr^2 E_u(1).
\end{equation}
\end{lem}

\begin{proof}
The second inequality follows at once from \eqref{3a2} (with $\lambda = 0$) for $u$. Next let
$v(x,t)=u(x,t)-\lambda_r(u)\, t$. Since $t$ is a solution for the constant coefficient
operator $L_0$, $v$ is a solution for $L_0$ as well in the domain in $T_1$,
 with $v(x,0)=0$ for all $x\in\Delta_1$. We claim that 
\begin{equation}\label{claim Dtv=0}
    \text{there exists some } (x',t')\in T_r \text{ for which } \Dt v(x',t')=0.
\end{equation} 
To see this, we observe first that $\Dt v(x,t)=\Dt u(x,t)-\fint_{T_r}\Dt u(x,t)dxdt$
has mean value $0$. 
Since $u$ is a solution of the constant-coefficient equation $L_0u=0$, 
$\Dt u$ is also a solution of the same equation. Therefore, by the De Giorgi-Nash-Moser theory, 
$\Dt u$ is continuous in $T_r$, and thus so is $\Dt v$. Then \eqref{claim Dtv=0} follows from the connectedness of $T_r$ and the mean value theorem.
Thanks to \eqref{claim Dtv=0}, $\underset{T_r}{\sup}\abs{\Dt v}\le\underset{T_r}{\osc}\Dt v$, and thus by \eqref{eq osc in use}
and because adding a constant does not change the oscillation, 
\begin{multline*}
    \fint_{T_r}\abs{\Dt v}^2\le \Bigl(\underset{T_r}{\osc}\Dt v\Bigr)^2=\Bigl(\underset{T_r}{\osc}(\Dt v+\lambda_r(u)-\lambda_1(u))\Bigr)^2\\
    =\Bigl(\underset{T_r}{\osc}\Dt(u-\lambda_1(u)\, t)\Bigr)^2
    \le C r^2\fint_{T_1}\abs{\nabla(u(x,t)-\lambda_1(u)\, t)}^2dxdt.
\end{multline*}
For the rest of the gradient, notice that
for $1\le j\le d$, 
$$\Dj v(x,t)=\Dj \br{v(x,t)+\lambda_r(u)\, t -\lambda_1(u)\, t},$$ 
and $\Dj v(x,0)=0$. Therefore,
\begin{multline*}
    \fint_{T_r}\abs{\Dj v}^2\le \Bigl(\underset{T_r}{\osc}\Dj \br{v(x,t)+\lambda_r(u)\, t -\lambda_1(u)\, t}\Bigr)^2\\
    \le C r^2\fint_{T_1}\abs{\nabla(u(x,t)-\lambda_1(u)\, t)}^2dxdt =Cr^2J_u(1).
\end{multline*}
Now \eqref{eq u0-lambda t} follows from the two estimates above.
\end{proof}

\begin{re}
The proof of Lemma~\ref{lem u0-lambda t} also works when we replace 
$J_u(r)$ in \eqref{eq u0-lambda t} with $\fint_{T_r}\abs{\nabla_{x,t}\br{u(x,t)-\lambda_s(u)\, t}}^2$, for any 
$0<s\le r$. That is, we also get that
\begin{equation}\label{3a8}
\fint_{T_r}\abs{\nabla_{x,t}\br{u(x,t)-\lambda_s(u)\, t}}^2dxdt 
\le C  r^2J_u(1).
\end{equation}
This may be a better estimate, since \eqref{3a2} says that for any $\lambda$, 
$$
J_u(r)
 \leq \fint_{T_r}\abs{\nabla\br{u(x,t)-\lambda\, t}}^2dxdt.
 $$
\end{re}

We will need a lower bound for the ratio 
$\frac{E_u(r)}{E_u(1)}$
for positive solutions of $L_0 u = 0$. 

\begin{lem}\label{lem lw bd}
Let the matrix $A_0$ be constant and satisfy 
the ellipticity condition \eqref{cond ellp},
set $L_0 = -\divg\br{A_0\nabla}$, and let $u$ be a positive solution to $L_0 u = 0$ in $T_1$ 
such that $u=0$ on $\Delta_1$. 
Then
\begin{equation}\label{eqlb}
E_u(r)
\ge C(1-C' r^2)E_u(1)
\qquad \text{ for } 
0<r<1/2,
\end{equation}
where $C$ and $C'$ are positive constants depending only on the dimension and $\mu_0$.
\end{lem}

Notice that when $r$ is small, the lower bound \eqref{eqlb} does not depend much on $r$. This
is better than what we would get by simply applying Lemma \ref{lem corkscrew} 
and the Harnack inequality to the positive solution $u$.  
The proof exploits the fact that $t$ is a solution for the constant-coefficient operator $L_0$ and the comparison principle.

\begin{proof}
Define $\lambda_0 =\Dt u(0,0)$.
Then by \eqref{eq osc in use}, 
\begin{equation*}
    \abs{\lambda_r(u)-\lambda_0}\le\underset{T_r}{\osc}\Dt u
    \le Cr\br{\fint_{T_1}\abs{\nabla u}^2}^{1/2}.
\end{equation*}
Since $t$ is a solution for $L_0$ that vanishes on $\Delta_1$, the comparison principle 
and Lemma \ref{lem corkscrew} give 
(with the corkscrew point $X_{x,t} = (x,t)$)
\[
\frac{u(x,t)}{t}
\geq C^{-1}u(X_{0,1})
\ge C^{-1}\br{\fint_{T_1}\abs{\nabla u}^2}^{1/2} 
\qquad \text{ for }  
(x,t)\in T_{1/2},
\]
which implies, 
by taking a limit and using the existence of $\nabla u$ at $0$, that
\[
\lambda_0=\Dt u(0,0)\ge C^{-1}\br{\fint_{T_1}\abs{\nabla u}^2}^{1/2}.
\]
Then
\begin{align*}
    E_u(r)\ge \lambda_r(u)^2
    \ge\frac{\lambda_0^2}{2}-(\lambda_r(u)-\lambda_0)^2
    \ge ((2C)^{-1}-C'r^2)\fint_{T_1}\abs{\nabla u}^2
\end{align*}
(use the fact that $a^2 \geq \frac{b^2}{2} - (a-b)^2$). 
This completes the proof of Lemma \ref{lem lw bd}.
\end{proof}

\subsection{Extension to general elliptic operators $L$}
We now return to a solution of our original equation  $Lu=0$, 
and compare it with solutions $u^0$ of $L_0 u^0=0$ of a constant coefficient
operator $L_0 = -\divg\br{A_0\nabla}$, with the same boundary data. 
For the moment we do not say who is the constant matrix $A_0$
(except that we require it to satisfy the ellipticity condition \eqref{cond ellp}),
but of course our estimates will be better if we choose a good approximation
of $A$ in $T_1$.

Even though it does not look like much, the next lemma is probably the central estimate of
this paper. We do not need $A_0$ to have constant coefficients here.

\begin{lem}\label{lem comp u u0}
Let $L = -\divg\br{A\nabla}$ and $L_0 = -\divg\br{A_0\nabla}$
be two elliptic operators, and assume that $A$ and $A_0$ satisfy the ellipticity condition \eqref{cond ellp}. Let $u$ be a solution to $Lu=0$ in $T_1$, with $u=0$ on $\Delta_1$,
and let $u^0$ be a solution of $L_0 u^0=0$ in $T_1$ with $u^0=u$ on $\bdy T_1$. 
Then there is some constant $C>0$ depending only on $d$ and the ellipticity constant $\mu_0$, such that 
\begin{equation} \label{eqmin}
\int_{T_1}\abs{\nabla u- \nabla u^0}^2 
\leq \mu_0^2 
    \min\set{\int_{T_1}\abs{A-A_0}^2\abs{\nabla u}^2dX,
    \int_{T_1}\abs{A-A_0}^2\abs{\nabla u^0}^2dX}.
\end{equation}
\end{lem}

\begin{proof} 
The solutions are in the space $W^{1,2}(T_1)$ by definition, and $u^0=u$ on the boundary should be interpreted as $u^0-u=0$ in the sense of $W^{1,2}(T_1)$, or equivalently, $u^0-u\in W_0^{1,2}(T_1)$.
So the existence of $u^0\in W^{1,2}(T_1)$ as above is guaranteed by the Lax-Milgram Theorem. 
Alternatively, it is possible to find $u^0$  because the trace of $u$ lies in $H^{1/2}(\d B)$. In addition, $u^0$ is nonnegative by the maximum principle.

Since $u-u^0$ lies in the set
$W^{1,2}_0$ of test functions allowed in Definition~\ref{def weak sol}, 

\begin{align*}
  \frac{1}{\mu_0} 
    \int_{T_1}\abs{\nabla(u-u^0)}^2&\le\int_{T_1} A\nabla(u-u^0)\cdot\nabla(u-u^0)=-\int_{T_1}A\nabla u^0\cdot\nabla(u-u^0)\\
    &=\int_{T_1}(A_0-A)\nabla u^0\cdot\nabla(u-u^0)\\
    &\le \frac{\mu_0}{2} 
    \int_{T_1}\abs{A-A_0}^2\abs{\nabla u^0}^2 +\frac{1}{2\mu_0}\int_{T_1}\abs{\nabla(u-u^0)}^2,
\end{align*}
where we use \eqref{cond ellp}, the fact that
$u$ is a solution of $\divg(A\nabla)u=0$ in $T_1$ (and $u-u^0$ vanishes on the boundary),
then the fact that $u^0$ is 
a solution of $\divg(A_0\nabla)u^0=0$ in $T_1$,
followed by the inequality $2ab \leq \mu_0 a^2 + \mu_0^{-1} b^2$. Then
\[
 \int_{T_1}\abs{\nabla(u-u^0)}^2\le \mu_0^2
 \int_{T_1}\abs{A-A_0}^2\abs{\nabla u^0}^2.
\]
This gives the bound by one of the expressions in the minimum in \eqref{eqmin}.
Interchanging the roles of $u$ and $u^0$, and $A$ and $A_0$, we also obtain the other bound.
\end{proof}

A similar proof also gives the following (which can be applied even if $A-A_0$ is not small).
\begin{lem}\label{lem u=u^0} 
Let $A$, $A_0$, $u$, and $u^0$ 
be as in Lemma \ref{lem comp u u0}. Then 
\begin{equation}\label{eq u=u^0}
\mu_0^{-4}\int_{T_1}\abs{\nabla u^0(X)}^2dX
\le\int_{T_1}\abs{\nabla u(X)}^2dX\le \mu_0^4\int_{T_1}\abs{\nabla u^0(X)}^2dX,
 \end{equation}
 where $\mu_0$ still denotes the ellipticity constant.
\end{lem}

We shall immediately see that $u$ being a solution is not necessary for the first inequality to hold, and similarly,  $u^0$ being a solution is not necessary for the second inequality. But the condition $u-u^0\in W_0^{1,2}(T_1)$ is essential.
\begin{proof} We estimate
\begin{align*}
    \mu_0^{-1}\int_{T_1}\abs{\nabla u}^2&\le\int_{T_1}A\nabla u\cdot\nabla u=\int_{T_1}A\nabla u \cdot\nabla(u-u^0)
    +\int_{T_1}A\nabla u \cdot\nabla u^0\\
    &=\int_{T_1}A\nabla u \cdot\nabla u^0\le\mu_0\br{\int_{T_1}\abs{\nabla u}^2}^{1/2}\br{\int_{T_1}\abs{\nabla u^0}^2}^{1/2}.
\end{align*}
Hence,
\[
  \int_{T_1}\abs{\nabla u}^2\le \mu_0^4\int_{T_1}\abs{\nabla u^0}^2.
\]
The left-hand side of \eqref{eq u=u^0} follows from the same argument, interchanging the roles of $u$ and $u^0$, $A$ and $A_0$, respectively. 
\end{proof}

Let us announce how we intend to estimate the right-hand side of \eqref{eqmin}.
The simplest would be to estimate $\abs{A-A_0}^2$ in $L^\infty$ norm and use the
$L^2$ norm of $\nabla u$, but if we do this we will get quantities that do not seem
to be controlled even by the $\alpha_\infty$ of \eqref{1a4}. So instead we decide to use
the quantity
\begin{equation} \label{3a15}
\gamma(x,r) = \inf_{A_0 \in \cA_0(\mu_0)}
\bigg\{\fint_{(y,s) \in T(x,r)} |A(y,s)-A_0|^2dyds \bigg\}^{1/2},
\end{equation}
where as before the infimum is taken over the class $\cA_0(\mu_0)$ of constant
matrices $A_0$ that satisfy the ellipticity condition \eqref{cond ellp}.
Notice that the domain of integration fits the domain of integration of \eqref{eqmin},
but it is larger than what we have in \eqref{1a5}. Nonetheless, the following lemma, 
to be proved in the next section,
will allow us to use the $\gamma(x,r)$.

\begin{lem}\label{l316}
If the matrix-valued function $A$ satisfies the weak DKP condition
of Definition \ref{d1a6}, with constant $\varepsilon > 0$,
then $\gamma(x,r)^2 \frac{dxdr}{r}$ is Carleson measure on $\R^{d+1}_+$, with norm
\begin{equation}\label{3a17}
\norm{\gamma(x,r)^2 \frac{dxdr}{r}}_{\mathcal{C}} 
\leq C \cN_2(A) \le C \epsilon,
\end{equation}
where $\cN_2(A) = \norm{\alpha_2(x,r)^2 \frac{dxdr}{r}}_{\mathcal{C}}$
as in \eqref{1a7}, and 
\begin{equation}\label{3a18}
\gamma(x,r)^2 \leq C \cN_2(A) \le C \epsilon 
\quad \text{ for } (x,r) \in \R^{d+1}_+.
\end{equation}
Here $C$ depends only on $d$ and $\mu_0$.
\end{lem}

See the next section for the proof. 

\ms 
Since we do not have a small $L^\infty$ control on $A$, we need a better estimate on $\nabla u$. This will be achieved  by reverse H\"older estimates (e.g. Lemma \ref{lem RH}), which gives us an exponent $p>2$ that depends only on $d$ and $\mu_0$. We first state the needed estimate for the unit box $T_1$.

\begin{lem}\label{l319}
Let $u$ be a positive solution to $Lu=0$ in $T_5$, with $u=0$ on $\Delta_5$, 
choose a constant matrix $A_0 \in \cA_0(\mu_0)$
that attains the infimum in the definition \eqref{3a15} of $\gamma(0,1)$,
and let $u^0$ be as in Lemma \ref{lem comp u u0} (with this choice of $A_0$). Then 
for any $\delta>0$,
\begin{equation} \label{3a20}
\int_{T_1} \abs{\nabla u- \nabla u^0}^2 dX
\leq \br{\delta+C_\delta \gamma(0,1)^2} E_u(1),
\end{equation}
where $C_\delta$ depends on $d$, $\mu_0$, and $\delta$. 
\end{lem}

\begin{proof}
We discussed the existence of $u^0$ when we proved Lemma \ref{lem comp u u0}.
We start from \eqref{eqmin}, which reads
\begin{equation} \label{3a21}
\int_{T_1} \abs{\nabla u- \nabla u^0}^2 
\leq C \int_{T_1}\abs{A-A_0}^2\abs{\nabla u}^2.
\end{equation}
Let us cut off and consider first the set 
\[Z:=\set{X\in T_1: \abs{\nabla u(X)}^2 \leq KE_u(1)},\]
with $K>0$ to be determined soon. We pull out the gradient and get a contribution
\begin{equation}\label{eq Z1}
   \int_Z\abs{A-A_0}^2\abs{\nabla u}^2 \leq KE_u(1)\int_{Z}\abs{A-A_0}^2
  \leq K \gamma(0,1)^2 E_u(1). 
\end{equation}
In the region $T_1\sm Z$ where $\abs{\nabla u}^2> K E_u(1)$, we see that
\[\abs{\nabla u}^2=\abs{\nabla u}^p \abs{\nabla u}^{2-p} 
\leq \abs{\nabla u}^p (K E_u(1))^{\frac{2-p}{2}},\]
where $p>2$ and will be chosen as in Lemma \ref{lem RH}. Then 
\begin{equation}\label{eq Z2}
  \int_{T_1\sm Z} \abs{A-A_0}^2\abs{\nabla u}^2
\leq 2\mu_0^2 \int_{T_1\sm Z} \abs{\nabla u}^2
\leq 2\mu_0^2 (K E_u(1))^{\frac{2-p}{2}} \int_{T_1} \abs{\nabla u}^pdX.
\end{equation}
We required $u$ to be a nice solution in the larger set $T_5$, so that we can use 
the following estimates from Section \ref{sec nota def}.
First, 
$$\big\{\fint_{T_1} \abs{\nabla u}^pdX\big\}^{\frac{2}{p}}
\leq C \fint_{T_2} \abs{\nabla u}^2dX
$$ 
by Lemma~\ref{lem RH}. Now we apply Lemma \ref{lem corkscrew} to $T_2$ 
(with $X_2 = (0,2)$) and later $T_1$ (with $X_1 = (0,1)$), to find that
\begin{equation*}
\fint_{T_2} \abs{\nabla u}^2 \leq C u^2(X_2) \leq C u^2(X_1) 
\leq C\fint_{T_1} \abs{\nabla u}^2,
\end{equation*}
where the intermediate inequality follows from Harnack's inequality. From these estimates and \eqref{eq Z2}, the contribution from $T_1\sm Z$ is 
\begin{equation*}
     \int_{T_1\sm Z} \abs{A-A_0}^2\abs{\nabla u}^2\le CK^{\frac{2-p}{2}}E_u(1).
\end{equation*}
Now we choose $K$ so that $CK^{\frac{2-p}{2}}=\delta$, and the desired estimate \eqref{3a20} follows at once. 
\end{proof}

We now have enough information to derive the same sort of decay estimates for 
the non-affine part of our solution $u$ that we proved, at the beginning
of this section, for solutions $u^0$ of constant coefficient operators.
We start with an analogue of Lemma \ref{lem u0-lambda t}.

\begin{lem}\label{lem Jur}
Let $u$ be a solution to $Lu=0$ in $T_1$ with $u=0$ on $\Delta_1$. Then 
for $0<r<1/4$,
\begin{equation} \label{Jur}
J_u(r) \leq C \br{r^2+K^{\frac{2-p}{2}}r^{-d-1}}J_u(1) 
+ \frac{C_K}{r^{d+1}}\gamma(0,1)^2E_u(1),
\end{equation}
where $K>0$ is arbitrary, $p=p(d,\mu_0)>2$, $C$ depends only on $d$, $\mu_0$ and $p$, and $C_K$ depends additionally on $K$.
\end{lem}

Notice that we do not require the positivity of $u$ yet, which is why we don't use Lemma \ref{l319} for the moment.

\begin{proof}
We write $u$ as affine plus orthogonal on $T_1$, i.e.
\[
u(x,t)=v(x,t)+\lambda_1(u)t.
\]
Note that $\lambda_1(u)^2\le E_u(1)$, and $E_v(1)=J_u(1)$.

Choose a constant matrix $A_0 \in \cA_0(\mu_0)$
that attains the infimum in the definition \eqref{3a15} of $\gamma(0,1)$, and let 
$L_0=-\divg{A_0\nabla}$ as usual. 
Now consider the $L_0$-harmonic extension to $T_{1/2}$ of the restriction of $u$ to $\bdy T_{1/2}$, which can be written as
\begin{equation}\label{ext u0}
 u_0(x,t) = v_0(x,t) + \lambda_1(u) t,   
\end{equation}
where we use the fact that $t$ is a solution of the constant-coefficient equation, and $v_0$ is the
$L_0$-harmonic extension of $v_{\vert \bdy T_{1/2}}$. These extensions are well-defined since $u$ is H\"older continuous on $\overline{T_{1/2}}$, and the Lax-Milgram Theorem guarantees the existence and uniqueness of the $W^{1,2}(T_{1/2})$ solution. In particular, $L_0 u_0=0$ in $T_{1/2}$, with $u_0=u$ on $\bdy T_{1/2}$. 

We claim that for any fixed $0<r<1/4$,
\begin{equation} \label{3a25}
J_u(r) \leq C r^2 J_u(1)
+ \frac{C}{r^{d+1}}\fint_{T_{1/2}}\abs{A(x,t)-A_0}^2\abs{\nabla u_0(x,t)}^2dxdt.
\end{equation}

To see this, we use the inequality $(a+b+c)^2 \leq 3(a^2 + b^2 + c^2)$ to write
\begin{multline}\label{est u-lambda r}
 J_u(r) = \fint_{T_r}\abs{\nabla\br{u-\lambda_r(u)\, t}}^2
    \le 3\fint_{T_r}\abs{\nabla(u_0 
    -\lambda_r(u_0)\, t)}^2
    \\+ 3\fint_{T_r}\abs{\nabla(u-u_0)}^2+3\fint_{T_r}\abs{\nabla(\lambda_r(u_0)\, t-\lambda_r(u)\, t)}^2, 
\end{multline}
where $\lambda_r(u_0) = \fint_{T_r} \partial_t u_0$ is defined as for $u$. Notice that
\begin{multline}\label{eq lambdas}
\fint_{T_r}\abs{\nabla(\lambda_r(u_0)\, t-\lambda_r(u)\, t)}^2
=(\lambda_r(u_0)-\lambda_r(u))^2 
=\br{\fint_{T_r}\br{\Dt u- \Dt u_0}dxdt}^2 \\ 
    \le\fint_{T_r}\abs{\nabla(u-u_0)}^2
    \le\frac{C}{r^{d+1}}\fint_{T_{1/2}}\abs{\nabla(u-u_0)}^2,
\end{multline}
simply enlarging the domain of integration. 
So by \eqref{est u-lambda r}, Lemma \ref{lem u0-lambda t} and Lemma~\ref{lem comp u u0}, 
\begin{multline}\label{eq u-lambda r cont}
J_u(r) \leq 3 \fint_{T_r}\abs{\nabla(u_0-\lambda_r(u_0)\, t)}^2 + \frac{C}{r^{d+1}}\fint_{T_{\frac{1}{2}}}\abs{\nabla(u-u^0)}^2\\
= 3 J_{u_0}(r) + \frac{C}{r^{d+1}}\fint_{T_{\frac{1}{2}}}\abs{\nabla(u-u_0)}^2
\leq C r^2 J_{u_0}(1/2) + \frac{C}{r^{d+1}}\fint_{T_{\frac{1}{2}}}\abs{\nabla(u-u_0)}^2\\
\leq C r^2 J_{u_0}(1/2) + \frac{C}{r^{d+1}}\fint_{T_{\frac{1}{2}}}
\abs{A-A_0}^2\abs{\nabla u_0}^2. 
 \end{multline}
However, the same sort of computation as above yields
\begin{multline*}
 J_{u_0}(1/2) = \fint_{T_{\frac{1}{2}}}\abs{\nabla(u_0-\lambda_{1/2}(u_0)t)}^2\\
    \le 3\fint_{T_{\frac{1}{2}}}\abs{\nabla(u-u_0)}^2+3\fint_{T_{\frac{1}{2}}}\abs{\nabla(u-\lambda_{1/2}(u)t)}^2+3(\lambda_{1/2}(u)-\lambda_{1/2}(u_0))^2\\
    \le C\fint_{T_{\frac{1}{2}}}\abs{\nabla(u-u_0)}^2+3\fint_{T_{\frac{1}{2}}}\abs{\nabla(u-\lambda_{1/2}(u)t)}^2\\
    = C\fint_{T_{\frac{1}{2}}}\abs{\nabla(u-u_0)}^2+3 J_u(1/2).
\end{multline*}
We plug this into \eqref{eq u-lambda r cont}, use the last part of \eqref{eq lambdas},
and get 
\begin{equation*} 
J_u(r) \leq C r^2 J_u(1/2)
+ \frac{C}{r^{d+1}}\fint_{T_{1/2}}\abs{A(x,t)-A_0}^2\abs{\nabla u_0(x,t)}^2dxdt.
\end{equation*}
Now the claim \eqref{3a25} follows because
\[
J_u(1/2)\le \fint_{T_{1/2}}\abs{\nabla(u(x,t)-\lambda_1(u)t)}^2dxdt\le CJ_u(1),
\]
where in the first inequality we have used that $\lambda_{1/2}(u)\,t$ is the best affine approximation in $T_{1/2}$ (see the discussion in Section \ref{subsec orth}).

Recall that $u_0$ is decomposed as in \eqref{ext u0}, and thus 
\begin{multline}\label{3b32}
     \fint_{T_{1/2}}\abs{A-A_0}^2\abs{\nabla u_0}^2\\
     \le 2\fint_{T_{1/2}}\abs{A-A_0}^2\abs{\nabla v_0}^2+2\lambda_1(u)^2\fint_{T_{1/2}}\abs{A-A_0}^2\abs{\nabla t}^2\\
    \le 2\fint_{T_{1/2}}\abs{A-A_0}^2\abs{\nabla v_0}^2+2E_u(1)\gamma(0,1)^2.
\end{multline}
We now estimate the first term on the right-hand side of \eqref{3b32}.
For $K>0$, consider the set
\[
Z_K:=\set{X\in T_{1/2}: \abs{\nabla v_0(X)}^2\le KE_u(1)}.
\]
The contribution of $Z_K$ to the integral is 
\[
\int_{Z_K}\abs{A-A_0}^2\abs{\nabla v_0}^2\le KE_u(1)\int_{Z_K}\abs{A-A_0}^2
\le CK\gamma(0,1)^2E_u(1).
\]
We are left with the complement of $Z_K$. 
As in \eqref{eq Z2} in the proof of Lemma \ref{l319}, we get that 
\begin{equation}\label{ZK cmpl}
     \int_{T_{1/2}\sm Z_K}\abs{A-A_0}^2\abs{\nabla v_0}^2\le C(KE_u(1))^{\frac{2-p}{2}}\int_{T_{1/2}}\abs{\nabla v_0}^p
\end{equation}
where $p>2$ will be chosen close to 2. To control the term $\int_{T_{1/2}}\abs{\nabla v_0}^p$, 
we use the following two reverse H\"older type estimates: for some $p=p(d,\mu_0)>2$ sufficiently close to $2$,
 \begin{align}
  \br{\int_{T_{1/2}}\abs{\nabla v_0}^p}^{1/p}
  &\lesssim \br{\int_{T_{1/2}}\abs{\nabla v_0}^2}^{1/2}
  +\br{\int_{T_{1/2}}\abs{\nabla v}^p}^{1/p},\label{eq RH1}\\
      \br{\int_{T_{1/2}}\abs{\nabla v}^p}^{1/p}&\lesssim \br{\int_{T_1}\abs{\nabla v}^2}^{1/2}+\abs{\lambda_1(u)}\br{\fint_{T_1}\abs{A-A_0}^p}^{1/p},\label{eq RH2}
 \end{align}
where the implicit constants depend on $d$, $\mu_0$ and $p$.
We postpone the proof of these two inequalities to the end of the proof of this lemma.

Now by \eqref{eq RH1} and \eqref{eq RH2}, we obtain
\[
\int_{T_{1/2}}\abs{\nabla v_0}^p\lesssim E_{v_0}(1/2)^{p/2}+E_{v}(1)^{p/2}+\abs{\lambda_1(u)}^p\fint_{T_1}\abs{A-A_0}^p.
\]
Since $v-v_0\in W_0^{1,2}(T_{1/2})$ and 
$v_0$ is $L_0$-harmonic, we have 
\[E_{v_0}(1/2)\le C_{\mu_0}E_v(1/2)\le CE_v(1)=CJ_u(1),\] 
where the first inequality comes from Lemma \ref{lem u=u^0}.
Notice also that
\[
\fint_{T_1}\abs{A-A_0}^p\le C_{\mu_0,p}\fint_{T_1}\abs{A-A_0}^2=C\gamma(0,1)^2.
\]
So our estimate on $\int_{T_{1/2}}\abs{\nabla v_0}^p$ can be simplified as 
\[
\int_{T_{1/2}}\abs{\nabla v_0}^p\lesssim J_{u}(1)^{p/2}+E_u(1)^{p/2}\gamma(0,1)^2.
\]
Plugging this into \eqref{ZK cmpl}, we get
\begin{multline*}
   \int_{T_{1/2}\sm Z_K}\abs{A-A_0}^2\abs{\nabla v_0}^2
   \le CK^{\frac{2-p}{2}}E_u(1)^{\frac{2-p}{2}}J_u(1)^{p/2}+CK^{\frac{2-p}{2}}\gamma(0,1)^2E_u(1)\\
   \le  CK^{\frac{2-p}{2}}J_u(1)+CK^{\frac{2-p}{2}}\gamma(0,1)^2E_u(1),
\end{multline*}
where in the last inequality we have used $E_u(1)\ge J_u(1)$, and thus $E_u(1)^{\frac{2-p}{2}}\le J_u(1)^{\frac{2-p}{2}}$. Combining this with the contribution on $Z_K$, we get
\[
 \int_{T_{1/2}}\abs{A-A_0}^2\abs{\nabla v_0}^2
   \le  CK^{\frac{2-p}{2}}J_u(1)+C\br{K+K^{\frac{2-p}{2}}}\gamma(0,1)^2E_u(1).
\]
From this and \eqref{3a25}, the desired estimate \eqref{Jur} follows.
\end{proof}

 We now prove the two H\"older type estimates. Let us first prove \eqref{eq RH2}.

{\em Proof of \eqref{eq RH2}.} 
Set $R_0=10^{-2}n^{-1/2}$ as before. 
For any $X_0=(x_0,t_0)\in T_{1/2}$, any $0<R\le R_0$, choose $\eta\in C_0^1(Q_{R}(X_0))$, with $\eta\equiv 1$ in $Q_{2R/3}(X_0)$, $\abs{\nabla \eta}\lesssim 1/R$. Here, $Q_R(X_0)$ is a cube centered at $X_0$ with side length $R$, and we shall write $Q_R$ for $Q_R(X_0)$ when this does not cause confusion.
Using $Lu=0$ in $T_1$, $v(x,t)=u(x,t)-\lambda_1(u)t$, and $L_0 t=0$, we have for any $w\in W_0^{1,2}(T_1)$,
\begin{multline}\label{RH2 eq1}
    0=\int_{T_1}A\nabla u\cdot\nabla w \,dxdt=\int_{T_1}A\nabla v\cdot\nabla w\,dxdt + \int_{T_1}A\nabla (\lambda t)\cdot\nabla w \,dxdt\\
    =\int_{T_1}A\nabla v\cdot\nabla w + \int_{T_1}(A-A_0)\nabla (\lambda t)\cdot\nabla w,
\end{multline}
where $\lambda=\lambda_1(u)$.

Now we choose 
$w(X)=v(X)\eta^2(X)$ when $t_0\le \frac{R}{2}$, and 
\(w=\Big(v-\fint_{Q_R}v(Y)dY\Big)\eta^2\) when $t_0>\frac{R}{2}$.
Notice that $v(x,0)=0$, and thus $w\in W_0^{1,2}(T_1)$ 
(because $Q_R \subset B_1$) as required. 
We plug $w$ into \eqref{RH2 eq1}, compute the derivatives, estimate some terms brutally,
and finally use Cauchy-Schwarz, and get the following estimates. 
\\
\textbf{Case 1:} $t_0\le \frac{R}{2}$.
Here we obtain 
\begin{multline*}
    \frac{1}{\mu_0}\int_{T_1}\abs{\nabla v}^2\eta^2dX\\
    \le\frac{1}{2\mu_0}\int_{T_1}\abs{\nabla v}^2\eta^2dX
    +C_{\mu_0}\int_{T_1}v^2\abs{\nabla \eta}^2dX
    + C_{\mu_0}\abs{\lambda}^2\int_{T_1}\abs{A-A_0}^2\eta^2dX.
\end{multline*}
Extending $v$ by zero below $t=0$, this yields 
\[
 \int_{Q_{2R/3}}\abs{\nabla v}^2dX
    \le\frac{C_{\mu_0}}{R^2}\int_{Q_R}v^2dX+ C_{\mu_0}\abs{\lambda}^2\int_{Q_R}\abs{A-A_0}^2dX.
\]
We apply the Poincar\'e-Sobolev inequality to control $\int_{Q_R}v^2dX$
and deduce from the above that
\begin{equation}\label{RH2 eq2}
    \fint_{Q_{2R/3}}\abs{\nabla v}^2dX
    \le C\br{\fint_{Q_R}\abs{\nabla v}^{\frac{2n}{n+2}}dX}^{\frac{n+2}{n}}+ C\abs{\lambda}^2\fint_{Q_R}\abs{A-A_0}^2dX.
\end{equation}
\textbf{Case 2:} $t_0> \frac{R}{2}$.
The same computation as in Case 1 gives
\begin{multline*}
    \int_{Q_{2R/3}}\abs{\nabla v}^2dX
    \le\frac{C}{R^2}\int_{Q_R}\Big|v(X)-\fint_{Q_R}v(Y)dY\Big|^2dX+ C\abs{\lambda}^2\int_{Q_R}\abs{A-A_0}^2dX.
\end{multline*}
Then by the Poinca\'re-Sobolev inequality,  
\eqref{RH2 eq2} holds again in this case. 

Now we apply \cite{giaquinta1983multiple} V. Proposition 1.1 to obtain
\[
\fint_{Q_{R_0/2}}\abs{\nabla v}^pdX
\le C\br{\fint_{Q_{R_0}}\abs{\nabla v}^2dX}^{\frac{p}{2}}
+ C\abs{\lambda}^p\fint_{Q_{R_0}}\abs{A-A_0}^pdX
\]
for some $p=p(d,\mu_0)>2$.

The desired estimate \eqref{eq RH2} follows as $T_{1/2}$ can be covered by finitely many $Q_{R_0/2}$.
\qed
\ms

 Now we turn to \eqref{eq RH1}.

{\em Proof of \eqref{eq RH1}.} 
We will use $L^p$ boundary estimates for solutions. Recall that $L_0 v_0=0$ in $T_{1/2}$, with $v_0-v\in W_0^{1,2}(T_{1/2})$. Set $R_0=10^{-2}n^{-1/2}$. Then by the boundary estimates in \cite{giaquinta1983multiple} p.154, we have for any $X_0\in T_{1/2}$,
\begin{multline*}
    \fint_{Q_{R_0/2}(X_0)\cap T_{1/2}}\abs{\nabla v_0}^p\lesssim \br{\fint_{Q_{R_0}(X_0)\cap T_{1/2}}\abs{\nabla v_0}^2}^{p/2}+\fint_{Q_{R_0}(X_0)\cap T_{1/2}}\abs{\nabla v}^p\\
    \lesssim \br{\fint_{T_{1/2}}\abs{\nabla v_0}^2}^{p/2}+\fint_{T_{1/2}}\abs{\nabla v}^p
\end{multline*}
for some $p>2$. Since $T_{1/2}$ can be covered by finitely many cubes $Q_{R_0/2}(X_0)$, we obtain \eqref{eq RH1}.
\qed
\ms

We now prove an analogue of Lemma \ref{lem lw bd} for positive solutions to $Lu=0$.

\begin{lem}\label{lem lwbd u}
Let $u$ be a positive solution of $Lu=-\divg(A\nabla)u=0$ 
in $T_5$, with $u=0$ on $\Delta_5$. 
Then for any $\delta>0$, $0<r<1/2$,
\begin{equation} \label{3a30}
E_u(r)
\ge\br{ \frac{1-C' r^2}{C}  
-\frac{C''\br{\delta+C_\delta\gamma(0,1)^2}}{r^{d+1}}}
E_u(1)
\end{equation}
where $C$, $C'$, $C''$ are positive constants depending only on $d$ and $\mu_0$.
\end{lem}

\begin{proof}
As before, we will only find this useful when the parenthesis is under control.  
Let $A_0$ and $u^0$ be as in Lemma \ref{l319}. By \eqref{3a20},
\begin{multline}\label{eq lwbd u est1}
\fint_{T_r}\abs{\nabla u}^2
\ge\frac{1}{2}\fint_{T_r}\abs{\nabla u^0}^2-\fint_{T_r}\abs{\nabla(u-u^0)}^2
\\ \ge\frac{1}{2}\fint_{T_r}\abs{\nabla u^0}^2
- \frac{1}{r^{d+1}}\fint_{T_1}\abs{\nabla(u-u^0)}^2
\\
\geq \frac{1}{2}\fint_{T_r}\abs{\nabla u^0}^2 
- \frac{C\br{\delta+C_\delta \gamma(0,1)^2}}{r^{d+1}}\fint_{T_1}\abs{\nabla u}^2.
\end{multline}
Divide both sides of \eqref{eq lwbd u est1} 
by $\fint_{T_1}\abs{\nabla u(X)}^2$,
and then observe that 
$$\fint_{T_1}\abs{\nabla u^0(X)}^2\approx\fint_{T_1}\abs{\nabla u(X)}^2 $$ 
by Lemma \ref{lem u=u^0};
this yields
\begin{multline*}
\frac{\fint_{T_r}\abs{\nabla u}^2}{\fint_{T_1}\abs{\nabla u}^2}
\geq \frac{1}{2}\frac{\fint_{T_r}\abs{\nabla u^0}^2}{\fint_{T_1}\abs{\nabla u}^2}
-\frac{C\br{\delta+C_\delta \gamma(0,1)^2}}{r^{d+1}}\\ 
\geq C^{-1}\frac{\fint_{T_r}\abs{\nabla u^0}^2}{\fint_{T_1}\abs{\nabla u^0}^2}
-\frac{C\br{\delta+C_\delta \gamma(0,1)^2}}{r^{d+1}}.
\end{multline*}
Since $u^0>0$ in $T_1$ (by the maximum principle), we can apply 
Lemma~\ref{lem lw bd} to $u^0$ and obtain the desired estimate.
\end{proof}

\ms
We are finally ready to prove the announced decay estimate for the quantity 
\begin{equation} \label{3a32}
\beta_u(x,r) = \frac{J_u(x,r)}{E_u(x,r)}
\end{equation}
(the proportion of non-affine energy) defined in \eqref{1a11}.
We just need to organize ourselves with the constants.

We intend to apply the estimates above, with a single value of $r = \tau_0$ 
which will be chosen small enough, depending on $d$
and $\mu_0$, and then we will 
require that
\begin{equation} \label{3a33}
\gamma(0,1) \leq \varepsilon_0,
\end{equation}
for some $\varepsilon_0 > 0$ that we shall choose momentarily, depending on 
$r = \tau_0$, $d$, and $\mu_0$.

Our first requirement for $r = \tau_0$ is that $C' r^2 < \frac12$ in \eqref{3a30}
(there will be another one of this type soon),
and we choose $\varepsilon_0$ and $\delta$ so small (depending on $\tau_0$)
that if \eqref{3a33} holds, then 
\[\frac{C''\br{\delta+C_\delta \gamma(0,1)^2}}{r^{d+1}} < \frac1{4C}\]
in \eqref{3a30}.
This way, \eqref{3a30} implies that
\begin{equation} \label{3a34}
E_u(r)
\ge  \frac{1}{4C}   E_u(1).
\end{equation}

Let $u$ be as in Lemma \ref{lem lwbd u}. 
We divide both sides of \eqref{Jur} by $E_u(r)$ and get that
\begin{equation} \label{3a36}
\beta_u(0,r) \leq C \br{r^2+K^{\frac{2-p}{2}}r^{-d-1}} \frac{J_u(1)}{E_u(r)} 
+\frac{C_K}{r^{d+1}}\gamma(0,1)^2\frac{E_u(1)}{E_u(r)}
\end{equation}
Then we choose $K$ to satisfy $K^{\frac{2-p}{2}}=r^{d+3}=\tau_0^{d+3}$, assume that \eqref{3a33} holds, apply \eqref{3a34}, and deduce from \eqref{3a36}
that (maybe with a larger constant $C$)
\begin{equation} \label{3a37}
\beta_u(0,\tau_0) \leq C \tau_0^2\beta_u(1) + C_{\tau_0}\gamma(0,1)^2.
\end{equation}
Finally we choose $\tau_0$ so small that (in addition to our earlier constraint) 
$C\tau_0^2 < \frac12$ in \eqref{3a37}, and finally choose $\varepsilon_0$ as above.

We recapitulate what we obtained so far in the next corollary.
Of course, by translation and dilation invariance, what was done with the unit box $T_1$ 
can also be done with any other
$T(x,R)$, $(x,R) \in \R^{d+1}_+$. We use the opportunity to state the general case,
which of course can easily be deduced from the case of $T_1$ by homogeneity
(or we could copy the proof).

\begin{cor}\label{cor itr}
We can find constants $\tau_0 \in (0,10^{-1})$ and $C > 0$
which depend only on $d$ and $\mu_0$, such that
if $u$ is a positive solution of $Lu=-\divg(A\nabla)u=0$ 
in $T(x,5R)$, with $u=0$ on $\Delta(x,5R)$, 
then 
\begin{equation} \label{3a39}
\beta_u(x,\tau_0 R) \leq  \frac12 \beta_u(x,R) + C \gamma(x,R)^2.
\end{equation}
\end{cor}

See \eqref{1a11} and \eqref{3a15} for the definitions of $\beta_u(x,\tau_0 R)$
and $\gamma(x,R)$.

\begin{proof}
The discussion above gives the result under the additional condition that 
$\gamma(x,R) \leq \varepsilon_0$. But we now have chosen $\tau_0$
and $\varepsilon_0$, and if $\gamma(x,R) > \varepsilon_0$, \eqref{3a39}
holds trivially (maybe with a larger constant), because $\beta_u(x,\tau_0 R) \leq 1$ by \eqref{3a3}.
\end{proof}

\begin{re}\label{r339}
As we remarked before, the complication of the decay estimate for $J_u(r)$ comes mainly from the lack of a small control of $\norm{A-A_0}_{L^\infty}$. 
If we knew $\gamma_\infty(x,R) \leq \varepsilon_1$, where 
\[\gamma_\infty(x,r)=\inf_{A_0\in\cA_0(\mu_0)}\sup_{T(x,r)}\abs{A-A_0},\]
then we 
could simplify the proof of Corollary \ref{cor itr} significantly. 

To see this, we start with an estimate similar to 
\eqref{3a25} \begin{equation}\label{3a25'} 
J_u(r) \leq C r^2 J_u(1)
+ \frac{C}{r^{d+1}}\fint_{T_{1/2}}\abs{A(x,t)-A_0}^2\abs{\nabla u(x,t)}^2dxdt,
\end{equation}
which can be obtained as \eqref{3a25}. Our estimate
for $\fint_{T_1}\abs{A-A_0}^2\abs{\nabla u}^2$ now becomes rather simple. We still choose 
$A_0$ as to minimize in the definition of $\gamma(0,1)$,
but observe that by Chebyshev, we can find $(x,t) \in T_1$
such that \[|A(x,t) - A_0| \leq C \gamma(0,1) \leq C\gamma_\infty(0,1).\]
Since $|A(y,s) - A(x,t)| \leq 2 \gamma_\infty(0,1)$ for $(y,s) \in T_1$,
we see that $|A-A_0| \leq C \gamma_\infty(0,1)\le C\eps_1$ on $T_1$. Then 
\begin{align*}
    \fint_{T_1}\abs{A-A_0}^2\abs{\nabla u}^2
 &\le 2\fint_{T_1}\abs{A-A_0}^2\abs{\nabla( u-\lambda_1(u)t)}^2+2\lambda_1(u)^2\fint_{T_1}\abs{A-A_0}^2\\
  &\le 2\fint_{T_1}\abs{A-A_0}^2\abs{\nabla( u-\lambda_1(u)t)}^2
    +2E_u(1)\fint_{T_1}\abs{A-A_0}^2 \\
 &\leq  2\varepsilon_1 J_u(1) + 2 \gamma(0,1)^2 E_u(1)
\end{align*}
and by \eqref{3a25'},
\begin{align*}
 J_u(r) \leq C \br{r^2 + \frac{C \varepsilon_1}{r^{d+1}}} J_u(1) + \frac{C\gamma(0,1)^2}{r^{d+1}} E_u(1).
\end{align*}
This is our analogue of \eqref{Jur}; the rest of the proof is the same.
\end{re}

\section{Carleson measure estimates}
\label{sec cm}

In this section we complete the proof of our two theorems. 
We already have our main decay estimate \eqref{3a39}, which says that $\beta_u(x,r)$
tends to get smaller and smaller, unless $\gamma(x,r)^2$ is large. This is a way of saying
that $\gamma^2$ dominates $\beta_u$, and it is not surprising that a Carleson measure estimate on the first function implies a similar estimate on the second one.
The fact that $\beta_u$ comes from a solution $u$ will not play any role in this argument.
See the second part of this section.

\subsection{Proof of Lemma \ref{l316}} \label{subs41}

Before we deal with decay, let us prove Lemma \ref{l316}, which is another fact about 
Carleson measures where $u$ plays no role.

Let $A$ be as in the statement. We want to show that 
$\gamma(x,r)^2 \frac{dxdr}{r}$ is Carleson measure on $\R^{d+1}_+$,
and our first move is to estimate $\gamma(x,r)$ in terms of the $\alpha_2(y,s)$.

For each pair $(x,r)$, we choose a constant matrix $A_{x,r}$ such that
\begin{equation} \label{4a1}
\fint_{W(x,r)} |A-A_{x,r}|^2 = \alpha_2(x,r)^2 .
\end{equation}
The interested reader may check that we can choose
the $A_{x,r}$ so that they depend on $(x,r)$ in a measurable way,
and in fact are constant on pieces of a measurable partition of $\R^{d+1}_+$,
maybe at the price of 
replacing $\alpha_2(x,r)^2$ in \eqref{4a1} with $2\alpha_2(x,r)^2$,
and making the $W(x,r)$ a little larger first to allow extra room to move $x$ and $r$.

Let $\Delta_0 = \Delta(x_0,r_0)$
be given; we want to estimate $\gamma(x_0,r_0)$,
and we try the constant matrix $A_0 = A_{x_0,r_0}$. Thus
\begin{equation} \label{4a2}
\gamma(x_0,r_0)^2 \leq \fint_{T_0} |A-A_0|^2 \leq C \fint_{Q_0} |A-A_0|^2,
\end{equation}
where we set $T_0 = T(x_0,r_0)$ and $Q_0 = \Delta(x_0,r_0) \times (0,r_0]$.
We will cut this integral into horizontal slices, using the radii
$r_m = \rho^m r_0$, $m \geq 0$. Let us choose $\rho = \frac45$,
rather close to $1$, to simplify the communication between slices.

We first estimate how fast the $A_{x,r}$ change. We claim that
\begin{equation} \label{4a3}
|A_{x,r} - A_{y,s}| \leq C \alpha_2(x,r) + C \alpha_2(y,s)
\quad \text{ when $|x-y| \leq \frac32 r$ 
and $\frac23 r \leq s \leq r$.}
\end{equation}
Indeed, with these constraints there is a box $R$ in $W(x,r) \cap W(y,s)$
such that $|R| \geq C^{-1} r^{d+1}$, and then 
\begin{align*} 
|A_{x,r} - A_{y,s}| 
&= \fint_R |A_{x,r} - A_{y,s}|  \leq \fint_R |A_{x,r} - A| + \fint_R |A - A_{y,s}| 
\\ &
\leq C \fint_{W(x,r)} |A_{x,r} - A| + C \fint_{W(y,s)} |A - A_{y,s}| 
\leq C \alpha_2(x,r) + C \alpha_2(y,s)
\end{align*}
 by the triangle inequality, the fact that $|R| \simeq |W(x,r)| \simeq |W(y,s)|$, and H\"older's inequality. 
 We can iterate this and get that for $y\in \R^d$ and $m \geq 0$,
 \begin{equation} \label{4a5}
|A_{y,r_m} - A_{y,r_0}| \leq C \sum_{j=0}^m \alpha_2(y,r_j).
\end{equation}
Now consider $y\in \Delta'_0 = \Delta(x_0,3r_0/2)$ and notice that by \eqref{4a3}, 
$|A_{y,r_0}-A_0| \leq C \alpha_2(y,r_0) + C \alpha_2(x_0,r_0)$, 
so \eqref{4a5} also yields
 \begin{equation} \label{4a6}
|A_{y,r_m} - A_0| \leq C \alpha_2(x_0,r_0) + C \sum_{j=0}^m \alpha_2(y,r_j).
\end{equation}
Set $H_m = \Delta_0 \times (r_{m+1},r_m]$ for $ m \geq 0$; thus $Q_0$ is the disjoint
union of the $H_m$. We claim that
\begin{equation} \label{4a7}
\int_{H_m} |A-A_0|^2 \leq C r_m  \alpha_2(x_0,r_0)^2  |\Delta_0|
+ C r_m \int_{\Delta'_0} \Big\{\sum_{j=0}^m \alpha_2(y,r_j) \Big\}^2 dy.
\end{equation}
We tried to discretize our estimates as late as possible, but this has to happen 
at some point. Cover $\Delta_0$ with disjoint cubes $R_i$ of sidelength 
$(10\sqrt d)^{-1}r_m$ that meet $\Delta_0$, and for each one choose a point $x_i \in R_i$ 
such that $\alpha_2(x_i,r_m)$ is minimal. Then set $A^i = A_{x_i, r_m}$ and 
$W_i = R_i \times (r_{m+1},r_m]$; notice that the $W_i$ cover $H_m$. 

The contribution of $R_i$ to the integral in \eqref{4a7} is
\begin{equation} \label{4a8}
\int_{W_i} |A(y,t)-A_0|^2dydt
\leq C \int_{W_i} |A(y,t)-A^i|^2 + |A^i-A_{y,r_m}|^2 + |A_{y,r_m}-A_0|^2 dydt.
\end{equation}
For the first term,
\begin{equation} \label{4a9}
\int_{W_i}|A(y,t)-A^i|^2 dydt \leq C |W(x_i,r_m)|  \alpha_2(x_i,r_m)^2 
\end{equation}
because $W_i \subset W(x_i,r_m)$ 
and by definition of $\alpha_2$. Next
\begin{align*}
\int_{W_i}|A^i-A_{y,r_m}|^2 dydt
\leq C \int_{W_i} (\alpha_2(x_i,r_m) + \alpha_2(y,r_m))^2 dydt
\leq C r_m \int_{R_i}  \alpha_2(y,r_m)^2 dy
\end{align*}
by \eqref{4a3} and because $\alpha_2(x_i,r_m)$, by the choice of $x_i$, is smaller. This integral is at least as large as
the previous one, again because $\alpha_2(x_i,r_m)$ is smaller.
When we sum all these terms over $i$, we get a contribution bounded by
$C r_m \int_{\Delta'_0} \alpha_2(y,r_m)^2$, which is dominated by the right
hand side of \eqref{4a7} (just keep the last term in the sum).
We are left with the third integral in \eqref{4a8}. But
$|A_{y,r_m}-A_0|$ is majorized in \eqref{4a6}, and the corresponding contribution,
when we sum over $i$, is also dominated by the right-hand side of \eqref{4a7}.
Our claim \eqref{4a7} follows. 

Because of \eqref{4a7} and the fact that the $H_m$ cover $Q_0$, we see that \eqref{4a2}
yields
\begin{equation} \label{4a11}
\gamma(x_0,r_0) \leq  C \fint_{Q_0} |A-A_0|^2 
\leq C |Q_0|^{-1} \sum_m \int_{H_m} |A-A_0|^2
\leq S_1 + S_2,
\end{equation}
where 
\begin{equation} \label{4a12}
S_1 = |Q_0|^{-1} \sum_m r_m  \alpha_2(x_0,r_0)^2  |\Delta_0| \leq C \alpha_2(x_0,r_0)^2,
\end{equation}
and 
\begin{equation} \label{4a13}
S_2 = |Q_0|^{-1} \sum_m r_m \int_{\Delta'_0} \Big\{\sum_{j=0}^m \alpha_2(y,r_j) \Big\}^2dy
\leq C \fint_{\Delta'_0} \sum_m \rho^m  \Big\{\sum_{j=0}^m \alpha_2(y,r_j) \Big\}^2dy
\end{equation}
because $r_m = \rho^m r_0$ and $|Q_0| \simeq r_0 |\Delta'_0|$.
We are about to apply Hardy's inequality, which says that for $1 < q < +\infty$,
\begin{equation} \label{4a14}
\sum_{m=0}^\infty \Big\{\frac{1}{m+1}\sum_{j=0}^m a_j \Big\}^q
\leq C_q \sum_{m} a_m^q
\end{equation}
for any infinite sequence $\{ a_m \}$ of nonnegative numbers. Here we take $q=2$ and 
$a_j = a_j(y) = \rho^{\frac{j}{4}}\alpha_2(y,r_j)$. Then
\begin{equation} \label{4a15}
\begin{aligned}
\sum_m \rho^m \Big\{\sum_{j=0}^m \alpha_2(y,r_j) \Big\}^2
&\leq \sum_m \rho^{m/2} \Big\{\sum_{j=0}^m \rho^{m/4} \alpha_2(y,r_j) \Big\}^2
\cr&\leq \sum_m \rho^{m/2} \Big\{\sum_{j=0}^m \rho^{j/4} \alpha_2(y,r_j) \Big\}^2
\cr& = \sum_m (m+1)^2 \rho^{m/2} \Big\{\frac{1}{m+1}\sum_{j=0}^m a_j \Big\}^2
 \leq C \sum_m a_m^2  
\end{aligned}
\end{equation}
so that
\begin{equation} \label{4a16}
S_2 \leq C  \fint_{\Delta'_0}\sum_m a_m^2(y) dy=C\sum_{m} \rho^{\frac{m}{2}} \fint_{\Delta'_0} \alpha_2(y,r_m)^2 dy.
\end{equation}
We return to \eqref{4a11}, use \eqref{4a12}, and see that
\begin{equation} \label{4a18}
\gamma(x_0,r_0)^2 \leq C \alpha_2(x_0,r_0)^2 
+ C \sum_{m} \rho^{\frac{m}{2}} \fint_{\Delta'_0} \alpha_2(y,\rho^m r_0)^2 dy
\end{equation}
We kept the squares because our Carleson measure condition is in terms of squares.
Recall that by assumption, $\alpha_2^2$ satisfies a Carleson measure condition, 
with norm $\cN_2(A)$. At this stage, deducing that the same thing holds for
$\gamma^2$ will only be a matter of applying the triangle inequality. We
write this because of the varying average in the second term of \eqref{4a18},
but not much will happen. Pick a surface ball $\Delta = \Delta(x_1,r_1)$. 
It is enough to bound 
\begin{equation} \label{4a19}
I = \int_\Delta \int_0^{r_1} \gamma(x,r)^2 \frac{dx dr}{r}
\leq C \int_\Delta \int_0^{r_1} \alpha_2(x,r)^2 \frac{dx dr}{r}
+ C \sum_m  \rho^{\frac{m}{2}} I_m,
\end{equation}
where 
\begin{equation} \label{4a20}
 I_m = \int_{x\in \Delta} \int_{r=0}^{r_1} \fint_{y \in \Delta(x,3r/2)} 
 \alpha_2(y,\rho^m r)^2 dy \frac{dx dr}{r}.
\end{equation}
Since 
\begin{equation} \label{4a21}
\int_\Delta \int_0^{r_1} \alpha_2(x,r)^2 \frac{dx dr}{r} \leq C \cN_2(A) r_1^d
\end{equation}
by definition, we may concentrate on $I_m$. Of course we apply Fubini. First
notice that $y \in \Delta' = \Delta(x_1,5r_1/2)$ when $y \in \Delta(x,3r/2)$ and $x\in \Delta$;
since $x\in \Delta(y,3r/2)$,
the integral in the dummy variable $x$ cancels with the normalization 
in the average, and we get that 
\begin{equation} \label{4a22}
 I_m = \int_{y\in \Delta'} \int_{r=0}^{r_1}  
 \alpha_2(y,\rho^m r)^2  \frac{dydr}{r}
 = \int_{y\in \Delta'} \int_{t=0}^{\rho^m r_1} \alpha_2(y,t)^2  \frac{dydt}{t},
\end{equation}
where the second identity is a change of variable 
(and we used the invariance of $\frac{dt}{t}$ under dilations). 
The definition also yields $I_m \leq C \cN_q(A) r_1^d$. So we can sum the series,
and we get that $I \leq C \cN_q(A) r_1^d$. This completes our proof of \eqref{3a17}.

We still need to check the second statement \eqref{3a18} (the pointwise estimate),
and this will follow from the fact that $\gamma$ is not expected to vary too much.
Indeed, we claim that
\begin{equation} \label{4a23}
\gamma(x,r) \leq C \gamma(y,s) 
\quad\text{whenever $|x-y| \leq r$ and $2r \leq s \leq 3r$.}
\end{equation}
This is simply because $T(x,r) \subset T(y,s)$, so if $A$ is well approximated by a constant coefficient matrix $A_0$ in $T(y,s)$, this is also true in $T(x,r)$. Now we square, average,
and get that
\begin{equation} \label{4a24}
\begin{aligned}
\gamma^2(x,r) &\leq C \fint_{y\in \Delta(x,r)}\fint_{s\in (2r,3r)} \gamma^2(y,s) dyds
\cr&\leq C r^{-d} \int_{y\in \Delta(x,r)}\int_{s\in (2r,3r)} \gamma^2(y,s) \frac{dyds}{s}
\leq C ||\gamma^2(y,s) \frac{dyds}{s}||_{\mathcal{C}} \leq C \cN_2(A).
\end{aligned}
\end{equation}
This completes our proof of Lemma \ref{l316}.
\ep

\begin{re} \label{r425}
There is also a local version of Lemma \ref{l316}, with the same proof.
It says that if $\alpha_2(x,r)^2 \frac{dxdr}{r}$ is Carleson measure relative to
some surface ball $3\Delta_0$ (see Definition \ref{d13}) , then 
$\gamma(x,r)^2 \frac{dxdr}{r}$ is Carleson measure on $T_{\Delta_0}$, 
with norm
\begin{equation}\label{4a26}
\norm{\gamma(x,r)^2 \frac{dxdr}{r}}_{\mathcal{C}(\Delta_0)} 
\leq C \norm{\alpha_2(x,r)^2 \frac{dxdr}{r}}_{\mathcal{C}(3\Delta_0)}.
\end{equation}
As usual, $C$ depends only on $d$.
For this the simplest is to observe that since we use nothing more than the estimate \eqref{4a18},
and for \eqref{3a17} we only care about $(x_0,r_0) \in T_{\Delta_0}$,
we may replace $\alpha_2(y,t)$ with $0$ when $(y,t) \notin T_{3\Delta_0}$.
Then the replaced function $\alpha_2$ satisfies a global square Carleson measure estimate
and we can conclude as above.

The fact that 
\begin{equation}\label{4a27}
\gamma(x,r)^2 \leq C \norm{\alpha_2(x,r)^2 \frac{dxdr}{r}}_{\mathcal{C}(3\Delta_0)}
\end{equation}
for $(x,r)\in T_{\Delta_0}$ can be proved as \eqref{3a18} above, using the fact that 
\eqref{4a26} also holds for a slightly larger ball $\frac{11}{10} \Delta_0$.
\end{re}

\subsection{Proof of Theorems \ref{mt1} and \ref{mt2}} \label{subs42}

We will just need to prove Theorem~\ref{mt2}, which is more general.
Let the matrix $A$ be as in the statement of both theorems. 

We recently completed our proof of Corollary \ref{cor itr}, which says that 
\begin{equation} \label{4a28}
\beta_u(x,\tau_0 r) \leq  \frac12 \beta_u(x,r) + C \gamma(x,r)^2
\end{equation}
whenever $u$ is a positive solution of $Lu=-\divg(A\nabla)u=0$ 
in $T(x,5r)$, with $u=0$ on $\Delta(x,5r)$. 

In the statement of our theorems, $u$ is assumed to be a 
positive solution of $Lu=0$ in $T(x_0,R)$, with $u=0$ on $\Delta(x_0,R)$, 
so \eqref{4a28} holds as soon as $\Delta(x,5r) \subset \Delta(x_0,R)$.
We pick such a pair $(x,r)$ and iterate \eqref{4a28}; this yields
\begin{equation} \label{4a29}
\beta_u(x,\tau_0^k r) \leq  
2^{-k} \beta_u(x,r) + C \sum_{j=0}^{k-1} 2^{-j} \gamma(x, \tau_0^{k-j-1}r)^2.
\end{equation}
Hence (writing $r$ in place of  $\tau_0^{-k} r$)
\begin{equation} \label{4a30}
\beta_u(x, r) \leq  
2^{-k} \beta_u(x,\tau_0^{-k} r) + C \sum_{j=0}^{k-1} 2^{-j} \gamma(x, \tau_0^{-j-1}r)^2
\end{equation}
as soon as $\Delta(x,5 \tau_0^{-k} r) \subset \Delta(x_0,R)$.

\ms
We want to prove the Carleson bound \eqref{1a15} on $\beta_u$ in $\Delta(x_0,\tau R)$,
so we give ourselves a surface ball $\Delta = \Delta(y,r) \subset \Delta(x_0,\tau R)$.
We want to show that 
\begin{equation} \label{4a31}
\int_{T_\Delta}\beta_u(x,s)\frac{dxds}{s} \leq C \tau^a r^d + C \cN r^d,
\end{equation}
where we set $\cN = \norm{\alpha_2(x,r)^2 \frac{dxdr}{r}}_{\mathcal{C}(\Delta(x_0,R))}$.

Let us first check that 
\begin{equation} \label{4a32}
\beta_u(x,s) \leq C \tau^a  + C \cN \quad\text{ when $x\in \Delta$ and $0<s \leq r$.}
\end{equation}
When $\tau \geq 10^{-1}$, this is true just because $(x,s) \in T(x_0,\tau R)$ and
\eqref{3a3} says that $\beta_u(x,s) \leq 1$. Otherwise, let $k$ be the largest integer
such that $\tau_0^{-k} r < 10^{-1} R$ (notice that $k \geq 0$); then
$\Delta(x,5 \tau_0^{-k} r) \subset \Delta(x_0,R)$, so \eqref{4a30} holds.
In addition, all the intermediate radii $\tau_0^{-j-1}r$ are also smaller
than $10^{-1} R$, so $\gamma(x, \tau_0^{-j-1}r)^2\leq C \cN$
by \eqref{3a18} or \eqref{4a27} in Remark \ref{r425}. Then \eqref{4a30}
says that $\beta_u(x,s) \leq  2^{-k} + C \cN$, and \eqref{4a32} follows,
with a constant $a$ that depends only on $\tau_0$
(which itself depends only on $d$ and $\mu_0$).  This is because our choice of $k$ gives 
$\tau_0^{k+1}\le 10r/R\le 10\tau$.

Call $I$ the integral in \eqref{4a31}, and write $I = \sum_{k= -1}^\infty I_k$, with
\begin{equation} \label{4a33}
I_k = \int_{T_\Delta} \1_{\tau_0^{k+2} r <  s \leq \tau_0^{k+1} r }(s) 
\beta_u(x,s)\frac{dxds}{s}.
\end{equation}
We single out $I_{-1}$ because we do not have enough room for the argument
below when $\tau$ is large, but anyway we just need to observe that 
\begin{equation} \label{4a34}
I_{-1} \leq C (\tau^a  +  \cN) |\Delta| \int_{\tau_0 r}^r \frac{ds}{s} 
\leq C (\tau^a  +  \cN) r^d
\end{equation}
by \eqref{4a32}, which is enough for \eqref{4a31}. 
We are left with $k \geq 0$ and 
\begin{equation} \label{4a35}
I_k \leq \int_{x\in \Delta} \int_{s = \tau_0^{k+2} r}^{\tau_0^{k+1} r } 
\beta_u(x,s)\frac{dxds}{s}.
\end{equation}
Because of our small precaution, we now have that for $(x,s)$ in the domain
of integration, $\tau_0^{-k}s \leq \tau_0 r \leq 10^{-1}r \leq 10^{-1} \tau R$ 
(because we took $\tau_0 \leq 10^{-1}$), so 
$\Delta(x,5\tau_0^{-k}s) \subset \Delta(x_0,R)$ and we can apply \eqref{4a30}. 
In addition, all the surface balls $5\Delta(x, \tau_0^{-j-1}s)$ that arise from \eqref{4a30}
are  contained in $\Delta(x_0,R)$, so we will be able to use 
Remark \ref{r425} to estimate them as in Lemma~\ref{l316}. 
Thus 
\begin{equation} \label{4a36}
\begin{aligned}
I_k &\le \int_{x\in \Delta} \int_{s = \tau_0^{k+2} r}^{\tau_0^{k+1} r } 
\big[2^{-k} \beta_u(x,\tau_0^{-k} s) 
+ C \sum_{j=0}^{k-1} 2^{-j} \gamma(x, \tau_0^{-j-1}s)^2\big] \frac{dxds}{s}
\cr& \leq C 2^{-k} (\tau^a  +  \cN) r^d 
+ C \sum_{j=0}^{k-1} 2^{-j}
\int_{\Delta} \int_{\tau_0^{k+2} r}^{\tau_0^{k+1} r } \gamma(x, \tau_0^{-j-1}s)^2\frac{dxds}{s}
\cr& = C 2^{-k} (\tau^a  +  \cN) r^d 
+ C \sum_{j=0}^{k-1} 2^{-j} \int_{\Delta} 
\int_{\tau_0^{k-j+1} r}^{\tau_0^{k-j} r } \gamma(x, t)^2 \frac{dxdt}{t}
\end{aligned}
\end{equation}
where we set $t = \tau_0^{-j-1}s$ and use the invariance of $\frac{ds}{s}$.

Set $\ell = k-j$, which runs between $1$ and $+\infty$. And for each value of 
$\ell \geq 0$, we have that $\sum_{k, j ; k-j = \ell} 2^{-j} \leq 2$. 
Hence when we sum over $k$, we get that
\begin{align*}
\sum_{k \geq 0} I_k 
&\leq C \sum_{k\ge 0} 2^{-k} (\tau^a  +  \cN) r^d  + C \sum_{\ell \geq 1} \int_{\Delta} 
\int_{\tau_0^{\ell+1} r}^{\tau_0^{\ell} r } \gamma(x, t)^2 \frac{dxdt}{t}
\\ & = C (\tau^a  +  \cN) r^d  + C  \int_{\Delta} 
\int_{0}^{\tau_0 r} \gamma(x, t)^2 \frac{dxdt}{t}
\leq C  (\tau^a  +  \cN) r^d,
\end{align*}
by Lemma~\ref{l316} or Remark \ref{r425}.
This completes our proof of \eqref{4a31}, and the theorems follow.

\section{Proof of Corollary \ref{cor main}}\label{sec cor}

Let us first prove a Caccioppoli type result for solutions on Whitney balls. Since it is an interior estimate, 
it holds on any domain $\Omega\subset\Real^{d+1}$. For $X\in\Omega$, denote by $\delta(X)$ 
the distance of $X$ to $\bdy\Omega$.
\begin{lem}\label{lem CcpType}
  Let $A$ be a $(d+1)\times (d+1)$ matrix of real-valued functions on $\Real^{d+1}$ 
satisfying the ellipticity condition \eqref{cond ellp}, and for some $C_0\in(0,\infty)$,
\begin{equation}\label{Adist}
    \abs{\nabla A(X)}\delta(X)\le C_0 \qquad\text{for any } X\in\Omega.
\end{equation}
 Let $X_0\in\Omega\subset\Real^{d+1}$ be given, and $r=\delta(X_0)$. Let $u\in W^{1,2}(B_{r}(X_0))$ be a solution of $Lu=-\divg(A\nabla u)=0$ in $B_{r}(X_0)$. Then for any $\lambda\in\Real$,
 \begin{multline}
     \int_{B_{r/4}(X_0)}\abs{\nabla^2u(X)}^2dX
     \le \frac{C}{r^2}\int_{B_{r/2}(X_0)}\abs{\nabla u(X)-\lambda\, \mathbf{e}_{d+1}}^2dX\\
     +C\lambda^2\int_{B_{r/2}(X_0)}\abs{\nabla A(X)}^2dX,
 \end{multline}
 where $C$ depends only on $d$, $\mu_0$ and $C_0$.
\end{lem}

\begin{proof}
By \eqref{Adist}, $\abs{\nabla A(X)}\le 8C_0/r$ for any $X\in B_{7r/8}(X_0)$, which means $A$ 
is Lipschitz in $B_{7r/8}(X_0)$. So from \cite{gilbarg2015elliptic} Theorem 8.8, it follows that 
$u\in W^{2,2}(B_{\frac{3}{4}r}(X_0))$. Let $\vp\in C_0^\infty(B_{r/2}(X_0))$, with $\vp=1$ on $B_{r/4}(X_0)$, $\norm{\nabla\vp}_{L^\infty}\le \frac{C}{r}$. Write ``$\d$" to denote a fixed generic derivative. 
Since  $u\in W^{2,2}(B_{\frac{3}{4}}(X_0))$, 
$\d (u-\lambda t)\vp^2\in W_0^{1,2}(B_{r/2}(X_0))$ for any $\lambda\in\Real$.
Therefore, there exists $\set{v_k}\subset C_0^\infty(B_{r/2}(X_0))$ such that $v_k$ 
converges to $\d (u-\lambda t)\vp^2$ in $W^{1,2}(B_{r/2}(X_0))$. 
Set $I=\int \abs{\nabla\d u(X)}^2\vp(X)^2 dX$. Observe that for any $\lambda\in\Real$, 
\[I=\int_{\Real^{d+1}}\abs{\nabla\d( u(x,t)-\lambda\, t)}^2\vp(x,t)^2 dxdt.\]
By ellipticity, we have 
\begin{multline*}
    I\le \mu_0\int_{\Real^{d+1}}A(x,t)\nabla\d( u(x,t)-\lambda\, t)\cdot\nabla\d( u(x,t)-\lambda\, t)\vp(x,t)^2 dxdt\\
    =\mu_0\int_{\Real^{d+1}}A\nabla\d( u-\lambda\, t)\cdot\nabla\br{\d( u-\lambda\, t)\vp^2}dxdt\\
    -2\mu_0\int_{\Real^{d+1}}A\nabla\d( u-\lambda\, t)\cdot\nabla\vp\, \d( u-\lambda\, t)\vp dxdt\\
    =:\mu_0 I_1-2\mu_0 I_2.
\end{multline*}
For $I_2$, we use 
Cauchy-Schwarz to get
\begin{multline*}
     \abs{I_2}\le\mu_0 I^{1/2}\br{\int_{\Real^{d+1}}\abs{\d( u-\lambda\, t)}^2\abs{\nabla\vp}^2dxdt}^{1/2}\\
    \le\frac{1}{8}I+\frac{C_{\mu_0}}{r^2}\int_{B_{r/2}(X_0)}\abs{\nabla(u-\lambda\,t)}^2dxdt.
\end{multline*}
For $I_1$, we use the sequence $\set{v_k}$ and write
\begin{multline*}
    I_1^k:=\int_{\Real^{d+1}}A\nabla\d( u-\lambda\, t)\cdot\nabla v_k dxdt\\
    =\int_{\Real^{d+1}}\d\br{A\nabla( u-\lambda\, t)\cdot\nabla v_k}dxdt-\int_{\Real^{d+1}}A\nabla( u-\lambda\, t)\cdot\nabla \d v_k dxdt \\
    -\int_{\Real^{d+1}}\d A(x,t) \nabla( u-\lambda\, t))\cdot\nabla v_k dxdt.   
\end{multline*}
Note that the first term on the right-hand side vanishes because it is a derivative of a $W^{1,2}(\Real^{d+1})$ compactly supported function. Moreover, since $Lu=0$  
and
$\partial v_k\in C_0^\infty(B_{r/2}(X_0))$ is a valid test function, we have
\begin{equation*}
     I_1^k=\lambda\int_{\Real^{d+1}}A\nabla t\cdot\nabla \d v_k dxdt 
    -\int_{\Real^{d+1}}\d A(x,t) \nabla( u-\lambda\, t))\cdot\nabla v_k dxdt.
\end{equation*}
Let $\mathbf{a}_{d+1}$ be the last column vector of $A$, then we have
\[
    \int_{\Real^{d+1}}A\nabla t\cdot\nabla \d v_k dxdt
    =\int_{\Real^{d+1}}\mathbf{a}_{d+1}\cdot\nabla \d v_k dxdt=-\int_{\Real^{d+1}}\dv\mathbf{a}_{d+1}\, \d v_k dxdt.  
\]
Hence,
\begin{multline*}
    \abs{I_1}=\abs{\lim_{k\to\infty}I_1^k}
    \le\abs{\lambda\int_{\Real^{d+1}}\dv\mathbf{a}_{d+1}\, \d (\d u\vp^2) dxdt}\\
    + \abs{\int_{\Real^{d+1}}\d A(x,t) \nabla( u-\lambda\, t))\cdot 
    \nabla(\d (u-\lambda t)\vp^2) \, 
         dxdt}=: I_{11}+I_{12}.
\end{multline*}
For $I_{11}$, we use  Cauchy-Schwarz, $\abs{\dv\mathbf{a}_{d+1}}\le (d+1)\abs{\nabla A}$, and Young's inequality to get
\begin{multline*}
    I_{11}\le\abs{\lambda}\int_{\Real^{d+1}}\abs{\dv\mathbf{a}_{d+1}}\d^2u\vp^2dxdt
    +2\abs{\lambda}\int_{\Real^{d+1}}\abs{\dv\mathbf{a}_{d+1}}\d(u-\lambda t)\vp\d\vp  dxdt\\
    \le\abs{\lambda}\br{\int\abs{\d^2u}^2\vp^2dxdt}^{1/2}\br{\int\abs{\dv\mathbf{a}_{d+1}}^2\vp^2dxdt}^{1/2}\\
    +2\abs{\lambda}\br{\int\abs{\d(u-\lambda t)}^2\abs{\nabla\vp}^2dxdt}^{1/2}\br{\int\abs{\dv\mathbf{a}_{d+1}}^2\vp^2dxdt}^{1/2}\\
    \le\frac{1}{8}I+\frac{C}{r^2}\int_{B_{r/2}(X_0)}\abs{\d(u-\lambda t)}^2dxdt+C\lambda^2\int_{B_{r/2}(X_0)}\abs{\nabla A}^2dxdt.
\end{multline*}
 For $I_{12}$, we have
 \begin{multline*}
     I_{12}\le\int_{\Real^{d+1}}
         |\d A(x,t)| \, |\nabla(u-\lambda t)| \, |\nabla(\d u)\vp^2| \, dxdt\\
      +2\int_{\Real^{d+1}}\abs{\d A(x,t)\nabla(u-\lambda t)\cdot\nabla\vp \d(u-\lambda t)\vp}dxdt\\
     \le I^{1/2}\br{\int_{B_{r/2}(X_0)}\abs{\d A}^2\abs{\nabla(u-\lambda t)}^2dxdt}^{1/2}\\
     +\frac{C}{r}\int_{B_{r/2}(X_0)}\abs{\d A}\abs{\nabla(u-\lambda t)}^2dxdt.
 \end{multline*}
By \eqref{Adist}, and 
because for any $X\in B_{r/2}(X_0)$, $\delta(X)\ge r/2$, one sees
\[
I_{12}\le \frac{1}{8}I+ \frac{C(d, C_0)}{r^2}\int_{B_{r/2}(X_0)}\abs{\nabla(u-\lambda t)}^2dxdt. 
\]
Collecting all the estimates, we can hide $I$ to the left-hand side and obtain the desired estimate.
\end{proof}

Let us point out that the assumption \eqref{Adist} on $A$ in Lemma \ref{lem CcpType} is harmless, as it is a consequence of the classical DKP condition \eqref{eq DKP}. We are now ready to prove Corollary \ref{cor main}.

\begin{proof}[Proof of Corollary \ref{cor main}]
Observe that \eqref{eq DKP} implies $\abs{\nabla A(x,t)}t\le CC_0$ for any $(x,t)\in \Real^{d+1}_+$ for some $C$ depending only on the dimension. 

Fix $\Delta\subset\Delta(x_0,R/2)$. Consider any $(x,3r)\in T_\Delta$, and write $X=(x,3r/2)$. Let  $\lambda_{x,3r}=\lambda_{x,3r}(u)$ be defined as in \eqref{1a9}. By Lemma \ref{lem CcpType},
\begin{multline*}
     \fint_{B_{r/4}(X)}\abs{\nabla^2u(y,t)}^2dydt
     \le \frac{C}{r^2}\fint_{B_{r/2}(X)}\abs{\nabla( u(y,t)-\lambda_{x,3r}t)}^2dydt\\
     +C\lambda_{x,3r}^2\fint_{B_{r/2}(X)}\abs{\nabla A(y,t)}^2dydt.
\end{multline*}
Notice that $B_{r/2}(X)\subset W(x,2r)=\Delta(x,2r)\times(r,2r]$ and $B_{r/2}(X)\subset T(x,3r)$. Hence we can enlarge the region of the integrals on the right-hand side and then multiply both sides by $u(x,3r)^{-2}r^3$ to get
\begin{multline*}
    \frac{\fint_{B_{r/4}(X)}\abs{\nabla^2u(y,t)}^2dydt}{ u(x,3r)^2}r^3
    \le  \frac{C r \fint_{T(x,3r)}\abs{\nabla( u(y,t)-{\lambda_{x,3r}}t)}^2dydt}{u(x,3r)^2}\\
     +\frac{Cr^3\lambda_{x,3r}^2}{u(x,3r)^2}\fint_{W(x,2r)}\abs{\nabla A(y,t)}^2dydt.
\end{multline*}
By Lemma \ref{lem corkscrew}, and 
then the definitions \eqref{1a8}-\eqref{1a10} of $\wt\alpha(x,r)$, $\lambda_{x,3r}$ and $\beta_u(x,3r)$,
\begin{multline*}
    \frac{\fint_{B_{r/4}(X)}\abs{\nabla^2u(y,t)}^2dydt}{ u(x,3r)^2} \; r^3\\ 
 \le \frac{C\fint_{T(x,3r)}\abs{\nabla( u(y,t)-\lambda_{ x,3r}t)}^2dydt}{r\fint_{T(x,3r)}\abs{\nabla u(y,t)}^2dydt}
    +\frac{C\br{\fint_{T(x,3r)}\d_t  u(y,t)dydt}^2 \wt\alpha(x,2r)^2}
    {r\fint_{T(x,3r)}\abs{\nabla u(y,t)}^2dydt}
      \\
     \le \frac{C\beta_u(x,3r)}{r}+\frac{C\wt\alpha(x,2r)^2}{r}.
\end{multline*}
Now we apply Theorem \ref{mt1} and the DKP assumption \eqref{eq DKP} and get
\begin{multline}\label{eq corpf1}
    \int_{T_\Delta}\frac{\fint_{B_{r/4}(X)}\abs{\nabla^2u(y,t)}^2dydt}{u(x,3r)^2} \; 
    r^3dxdr\\
    {\le C\int_{T_\Delta}\beta_u(x,3r)\frac{dxdr}{r}+C\int_{T_{\Delta}}\wt\alpha(x,2r)^2\frac{dxdr}{r}}
    \le C_{d,\mu_0}(1+C_0)\abs{\Delta}.
\end{multline}
We now use Fubini and
 Harnack's inequality to obtain a lower bound for the left-hand side of \eqref{eq corpf1}. 
By Fubini,
\begin{multline*}
    \int_{T_\Delta}\frac{\fint_{B_{r/4}(X)}\abs{\nabla^2u(y,t)}^2dydt}{u(x,3r)^2} \; 
    r^3dxdr\\
    = C_d\int_{(y,t)\in \Real^{d+1}_+}\abs{\nabla^2u(y,t)}^2\int_{(x,r)\in T_\Delta}\1_{B_{r/4}(X)}(y,t)\frac{r^{2-d}}{u(x,3r)^2} \; 
    dxdrdydt.
\end{multline*}
Observe that if $\abs{(y,t)-(x,3r/2)}\le \frac{t}{7}$, then $t\approx r$, $t\le \frac{7r}{4}$, and the latter implies that $\1_{B_{r/4}(X)}(y,t)=1$. 
So the right-hand side is bounded from below by
\[
    C_d\int_{(y,t)\in T_\Delta}\abs{\nabla^2u(y,t)}^2
 \int_{(x,r) ; (x,3r/2)\in B_{t/7}(y,t)}\frac{r^{2-d}}{u(x,3r)^2} \; dxdrdydt.
   \]
By Harnack, $u(x,3r)\le C u(y,t)$ 
when $(x,3r/2)\in B_{t/7}(y,t)$. 
Hence
\begin{multline*}
    \int_{T_\Delta}\frac{\fint_{B_{r/4}(X)}\abs{\nabla^2u(y,t)}^2dydt}{u(x,3r)^2}r^3dxdr
    \ge C_d\int_{(y,t)\in T_\Delta}\abs{\nabla^2u(y,t)}^2\frac{t^3}{u(y,t)^2}dydt.
\end{multline*}
From this and \eqref{eq corpf1}, the desired result follows. 
\end{proof}

\begin{re}
If we apply the more precise estimate \eqref{1a15} in \eqref{eq corpf1}, we can get the following stronger result. 
For $\tau\in(0,1/2)$, we have 
\[
\norm{\frac{\abs{\nabla^2 u(x,t)}^2t^3}{u(x,t)^2}dxdt}_{\C(\Delta(x_0,\tau R)}\le C\tau^a+C\norm{\wt\alpha(x,r)^2\frac{dxdr}{r}}_{\C(\Delta(x_0,R))},
\]
for some $C$ and $a>0$ depending only on $d$ and $\mu_0$. As a consequence, if $u$ is the Green function with pole at infinity (see Lemma \ref{lem emGreen_infinity} for the definition), then we have that 
\[
\norm{\frac{\abs{\nabla^2 G^\infty(x,t)}^2t^3}{G^\infty(x,t)^2}dxdt}_{\C}\le C\norm{\wt\alpha(x,r)^2\frac{dxdr}{r}}_{\C}.
\]
\end{re}

\section{Optimality}\label{sec optm}
In this section, 
we construct an operator that does not satisfy the DKP condition and such that  $\beta_{G^\infty}(x,r)\frac{dxdr}{r}$ 
fails to be a Carleson measure. Moreover, we find a sequence of operators $\{L_n\}$ that satisfy the DKP condition with constants increasing to infinity as $n$ goes to infinity, and for any fixed $1<R_0<\infty$,
$\norm{\beta_n(x,r)\frac{dxdr}{r}}_{\mathcal{C}(\Delta_{R_0})}\ge C(n-1)$, where 
$\beta_n(x,r)=\beta_{G_n^\infty}(x,r)$, and $G_n^\infty$ is the Green function with pole at infinity for $L_n$. 
A similar construction is used in \cite{DM2020} Remark 3.2 and \cite{DFM2020}.
As we shall see, it is very simple to get a bad oscillating behaviour for $G^\infty$ in the vertical
direction; it is typically harder to get oscillation in the horizontal variables, as would be needed for bad
harmonic measure estimates. 

 Let us give the precise definition of the Green function with pole at infinity. One can prove the following lemma as in \cite{kenig1999free}, 
Lemma 3.7.
\begin{lem}\label{lem emGreen_infinity}
Let $L=-\divg{A\nabla}$ be an elliptic operator on $\Real^{d+1}_+$. Then there exists a unique function $U\in C(\overline{\Real^{d+1}_+})$ such that 
\[
\begin{cases}
L^\top U = 0 \quad \text{in } \Real^{d+1}_+\\
U>0 \quad \text{in }\Real^{d+1}_+\\
U(x,0)=0 \quad \text{for all }x\in\Real^d,\\
\end{cases}
\]
and $U(0,1)=1$. 
We call the unique function $U$ {\it the Green function with pole at infinity for $L$}.
\end{lem}

Let $A(x,t)=a(t)I$ for $(x,t)\in\Real^{d+1}_+$, where $I$ is the $d+1$ identity matrix, and $a(t)$ is a positive scalar function on $\R_+$. Let $L=-\divg{A(x,t)\nabla}$. We claim that the Green function with pole at infinity 
for $L$ in $\Real^{d+1}_+$ is 
(modulo a harmless multiplicative constant)
\begin{equation}\label{Ginfty}
 G(x,t)=g(t) \qquad\text{with}\quad g(0)=0,\quad g'(t)=\frac{1}{a(t)}.   
\end{equation}
In fact, it is easy to check that $L^\top G=0$ in $\Real^{d+1}_+$, $G(x,0)\equiv 0$, 
and the uniqueness of $G^\infty$ does the rest. 
The derivatives of $G$ are simple. They are
\begin{equation}\label{de G}
    \nabla_xG(x,t)=0, \qquad \Dt G(x,t)=\frac{1}{a(t)}.
\end{equation}
Now we set 
\[
a(t)=
\begin{cases}
\frac{3}{2} \qquad\text{when }t\ge 2^{100},\\
1 \qquad\text{when }2^{2k}+c_02^{2k-1}\le t\le 2^{2k+1}-c_02^{2k},\\
2 \qquad\text{when }2^{2k+1}+c_02^{2k}\le t\le 2^{2k+2}-c_02^{2k+1},
\end{cases}
\]
for all $k\in\mathbb{Z}$ with $k\le 49$, and $a(t)$ is smooth 
in the remaining strips $S_k = (2^k-c_02^{k-1},2^k+c_02^{k-1})$, 
with \[\abs{a'(t)}\le\frac{100}{c_02^k}\quad\text{ for } t\in S_k= (2^k-c_02^{k-1},2^k+c_02^{k-1}).\] Here, $c_0>0$ is a constant that will be taken sufficiently small and fixed. 
Additionally, we can make sure that $a(t)=\frac{3}{2}$ in a small neighborhood of $t=2^k$
to simplify our computations.

We construct the approximation of $a(t)$ as follows. Set
\[
a_n(t)=\begin{cases}
a(t)\qquad\text{when } t\ge 2^{-2n},\\
\frac{3}{2} \qquad\text{when }0<t<2^{-2n}.
\end{cases}
\]
Then $a_n$ converges to $a$ pointwise in $\Real^{d+1}_+$.

Let $L_n=-\divg{A_n(x,t)\nabla}=-\divg\br{a_n(t)\nabla}$, and let $G_n$ be the Green function with pole at infinity for $L_n$, whose formula are given in \eqref{Ginfty}.

We now compute the DKP constant for $A_n$. Notice that $\abs{\nabla A_n}\neq 0$ only in the strips near $2^k$ with width $c_02^k$ for $-2n\le k\le 100$, so it is easy to get the following estimate.
\begin{multline*}
    \norm{\sup_{(y,t)\in W(x,r)}\abs{\nabla A_n(y,t)}^2r dxdr}_{\mathcal{C}}\approx\norm{\abs{ a'_n(t)}^2tdxdt}_{\mathcal{C}}\\\approx\sum_{k=-2n}^{100}\frac{2^k}{(c_02^k)^2}c_02^k\approx\frac{2n+100}{c_0}.
\end{multline*}
Similarly, we can compute the DKP constant for $A$.
\begin{multline*}
    \norm{\sup_{(y,t)\in W(x,r)}\abs{\nabla A(y,t)}^2r dxdr}_{\mathcal{C}}\approx\norm{\abs{ a'(t)}^2tdxdt}_{\mathcal{C}}\approx c_0^{-1}\sum_{k=-\infty}^{100}1=\infty.
\end{multline*}
Now we turn to $\beta_n$. Recall the definition of $\beta(x,r)$ \eqref{1a11} and the simple expressions for the derivatives of $G_n$ \eqref{de G}.
Set $b_n(t)=\frac{1}{a_n(t)}$ and compute $\beta_n(x,r)$ with $T(x,r)$ replaced by $\Delta(x,r)\times(0,r)$
 in the definition of $\beta(x,r)$; then 
\begin{multline}\label{beta_n}
    \beta_n(x,r)=\frac{\int_{y\in \Delta(x,r)}\int_{t=0}^r\abs{\Dt G_n(y,t)-\fiint_{\Delta(x,r)\times(0,r)}\Dt G_n(y',t')dy'dt'}^2dtdy}
    {\int_{y\in \Delta(x,r)}\int_{t=0}^r\abs{\nabla G_n(y,t)}^2dtdy}\\
    =\frac{\int_0^r\abs{b_n(t)-\fint_0^rb_n(s)ds}^2dt}{\int_0^r\abs{b_n(t)}^2dt}. 
\end{multline}
The estimates with our initial definition of $T(x,r)$ would be very similar, or could be deduced from the
estimates with $\Delta(x,r)\times(0,r)$ because $T(x,r/10) \subset \Delta(x,r)\times(0,r) \subset T(x,10)$.

Notice that $\beta_n(x,r)=0$ when $r<2^{-2n}$.
We estimate $\norm{\beta_n(x,r)\frac{dxdr}{r}}_{\mathcal{C}(\Delta_{R_0})}$ for some fixed $R_0\ge 1$. For simplicity, we only do the calculation when $R_0<2^{100}$. 

The main observation is that for any $2^{-2n+2}\le r \leq R_0$, 
\begin{equation}\label{bn lwbd}
    \abs{b_n(t)-\fint_0^rb_n(s)ds}^2\ge\frac{1}{1000} \qquad\text{for }
      t\in [2^{-2n},r]\setminus (\cup_k S_k).
    \end{equation}
Once we have \eqref{bn lwbd}, we can obtain the lower bound for $\norm{\beta_n(x,r)\frac{dxdr}{r}}_{\mathcal{C}(\Delta_{R_0})}$ as follows. First, observe that the total measure of those $S_k$ that intersects $[2^{-2n},r]$ is controlled. Namely,
\[
\abs{\cup_k S_k\cap[2^{-2n},r]}\le\sum_{k=-2n}^{-2n+j+1}c_02^k\le c_02^{-2n+j+2}\le 4c_0r,
\]
where $j$ is the integer that $2^{-2n+j}\le r<2^{-2n+j+1}$. Therefore,
\[
\int_0^r\abs{b_n(t)-\fint_0^rb_n(s)ds}^2dt\ge\int_{[2^{-2n},r]\setminus \br{\cup_k S_k}}\frac{1}{1000}dt\ge\frac{\frac3{4}-4c_0}{1000}r=:C_0 r.
\]
On the other hand, we have 
\(
\int_0^r\abs{b_n(t)}^2dt\le r
\)
since $\abs{b_n}\le 1$. Then by the formula \eqref{beta_n} for $\beta_n$, we obtain
\[
\beta_n(x,r)\ge C_0 \qquad\text{for }r\in [2^{-2n+2}, R_0].
\]
So 
\begin{multline*}
    \sup_{0<R\le R_0}\frac{1}{R^d}\int_{\Delta_R}\int_0^R\beta_n(x,r)\frac{dxdr}{r}\ge\frac{\abs{\Delta_{R_0}}}{R_0^d}\int_{2^{-2n+2}}^{R_0}C_0\frac{dr}{r}\\
    =C_{d,c_0}\br{(2n-2)\ln 2+\ln R_0}\ge C_{d,c_0}(2n-2).
\end{multline*}

Now we justify \eqref{bn lwbd}. This is true simply because the average $\fint_0^rb_n(s)ds$ takes value strictly between $1$ and $\frac{1}{2}$, so when $t$ is away from the strips $S_k$, $b_n(t)$ should be different than $\fint_0^rb_n(s)ds$. We just need to make sure that the lower bound does not depend on $n$ in a way that would cancel the blow up.  

We first simplify our computation of $\fint_0^rb_n(s)ds$ by observing that we can take $c_0=0$. This is because if $c_0\neq0$, we can always require the average of $b_n$ in $(0,r)$ to be the same as the case when $b_n$ is not smoothed out (i.e. $c_0=0$), as long as $r$ does not lie in any strip $S_k$, by choosing our $a_n$ carefully. But if $r\in S_k$, this should not affect $\fint_0^rb_n(s)ds$ much if we take $c_0$ to be sufficiently small. 

Fix $2^{-2n+2}\le r\le R_0$. If $2^{2k_0}\le r<2^{2k_0+1}$ for some $k_0\in\mathbb{Z}$, then a direct computation shows
\begin{equation*}
    \fint_0^rb_n(s)ds=1+\frac{2^{-2n}}{2r}-\frac{2^{2k_0}}{3r}.
\end{equation*}
If $2^{2k_0}\le r<2^{2k_0+1}$ for some $k_0\in\mathbb{Z}$, then
\begin{equation*}
     \fint_0^rb_n(s)ds=\frac{1}{2}+\frac{2^{-2n}}{2r}+\frac{2^{2k_0+1}}{3r}.
\end{equation*}
Since $b_n$ is either $1$ or $1/2$ in $(0,r)\setminus S_k$, a case-by-case computation shows that for any $2^{-2n+2}\le r\le R_0$,
$\abs{b_n(t)-\fint_0^rb_n(s)ds}\ge\frac{1}{12}$ for $t\in [2^{-2n},r]\setminus S_k$. Then with $c_0>0$ sufficiently small, we have \eqref{bn lwbd}.

\bibliographystyle{alpha}
\bibliography{references}

\Addresses
\end{document}